\documentclass[11pt, oneside]{amsart}   	
\usepackage{geometry}                		
\usepackage{graphicx}				
\usepackage{amssymb}
\usepackage{epstopdf}
\usepackage{amsmath}
\usepackage{amsfonts}
\usepackage{fancyhdr} 
\usepackage{mathrsfs}
\usepackage{wrapfig}
\usepackage{bbold}
\usepackage{amsthm} 
\usepackage{enumitem}
\usepackage[toc]{appendix}
\usepackage{xcolor}

\usepackage[colorlinks=true, pdfstartview=FitV, linkcolor=blue, 
        citecolor=blue, urlcolor=blue]{hyperref}

\newtheorem{thm}{Theorem}
\newtheorem{prop}{Proposition}[section]
\newtheorem{lem}[prop]{Lemma}
\newtheorem{defn}[prop]{Definition}
\newtheorem{cor}[prop]{Corollary}

\newenvironment{pf}{\noindent {\bf\em Proof}.\ \ }{\hspace*{\fill}\rule{1.4ex}{1.4ex}\,}

\newcommand{\ddd}{\mathbb{D}}
\newcommand{\real}{\mathbb{R}}
\newcommand{\complex}{\mathbb{C}}
\newcommand{\integer}{\mathbb{Z}}

\newcommand{\torus}{\mathbb{T}}
\newcommand{\pro}{\mathfrak{p}}
\newcommand{\rat}{\mathbb{Q}}

\DeclareMathOperator{\diag}{diag}

\DeclareMathOperator{\supp}{supp}

\DeclareMathOperator{\spn}{span}
\DeclareMathOperator{\aut}{Aut}

\usepackage{scalerel,stackengine}
\stackMath
\newcommand\reallywidecheck[1]{%
\savestack{\tmpbox}{\stretchto{%
  \scaleto{%
    \scalerel*[\widthof{\ensuremath{#1}}]{\negmedspace\kern0pt\bigwedge\kern0pt}%
    {\rule[-\textheight/2]{1ex}{\textheight}}
  }{.5\textheight}%
}{0.5ex}}%
\stackon[1pt]{#1}{\scalebox{-1}{\tmpbox}}%
}

\title[Theorems of Szeg\H{o}-Verblunsky type]{Theorems of Szeg\H{o}-Verblunsky type in the multivariate and almost periodic settings}
\author{Peter C.~Gibson}

\date{February 20, 2022}							

\begin{document}

\begin{abstract}
The classical Szeg\H{o}-Verblunsky theorem relates integrability of the logarithm of the absolutely continuous part of a probability measure on the circle to square summability of the sequence of recurrence coefficients for the orthogonal polynomials determined by the measure.  The present paper constructs orthogonal polynomials on the torus of arbitrary finite dimension in order to prove theorems of Szeg\H{o}-Verblunsky type in the multivariate and almost periodic settings. The results are applied to the one-dimensional Schr\"odinger equation in impedance form to yield a new trace formula valid for piecewise constant impedance, a case where the classical trace formula breaks down.  As a byproduct, the analysis gives an explicit formula for the Taylor coefficients of a bounded holomorphic function on the open disk in terms of its continued fraction expansion. 
 \end{abstract}

\maketitle

\begin{center}
MSC 42C05 (32A10, 42B30, 42A75, 34L25)
\end{center}

\tableofcontents

\section{Introduction\label{sec-introduction}}

The present paper constructs orthogonal polynomials on the $d$-dimensional torus in order to extend Szeg\H{o}'s classical theorem concerning measures on the circle to higher dimensions, and to almost periodic functions. The underlying objective is a new trace formula for the one-dimensional Schr\"odinger equation with singular potential. 

An exotic Riemannian structure on the unit disk, called the scattering metric, is introduced as a key technical ingredient, facilitating representation of the Taylor coefficients of an arbitrary bounded holomorphic function on the open unit disk in terms of its continued fraction parameters---thus solving in explicit form a problem raised by Schur in 1917.

\subsection{Notation, objectives and background}
Write $\ddd$ for the open unit disk in the complex plane, $\overline{\ddd}$ for its closure, and $\torus$ for its boundary, the unit circle. For $d\geq1$, let $\mathcal{H}_d$ denote the set of all holomorphic functions 
\begin{equation}\label{Hd}
h:\ddd^d\rightarrow\overline{\ddd},
\end{equation}
endowed with the topology of uniform convergence on compact sets.  Any $h\in\mathcal{H}_d$ has an associated diagonal, $\diag h\in\mathcal{H}_1$, defined by the equation
\begin{equation}\label{diagonal}
\diag h(z)=h(z,\ldots,z)\qquad(z\in\ddd). 
\end{equation}
For every $h\in\mathcal{H}_d$, the multivariate version of Fatou's theorem \cite[Thms.~2.1.3(e), 2.3.1]{Ru:1969} guarantees the existence almost everywhere on $\torus^d$ of radial limits
\[
h^\circ(z):=\lim_{s\nearrow 1}h(sz)\qquad (z\in\torus^d). 
\]
Let $\tau$ denote normalized Lebesgue measure on $\torus^d$, so that, for $f\in L^1(\torus^d)$, 
\begin{equation}\label{tau-notation-0}
\int_{\torus^d}f\,d\tau=\frac{1}{(2\pi)^d}\displaystyle\rule{-27pt}{0pt}\int\limits_{(e^{i\theta_1},\ldots,e^{i\theta_d})\in\torus^d}\rule{-27pt}{0pt}f(e^{i\theta_1},\ldots,e^{i\theta_d})\,d\theta_1\cdots d\theta_d.
\end{equation} 
The same symbol $\tau$ will be used for all dimensions, with the pertinent value of $d$ being evident from context. 
In the present paper sequences are represented as functions on the non-negative integers $\integer_+$, or positive integers $\integer_{>0}$, with function values being indicated by a subscript.  Thus for example, $r:\integer_+\rightarrow\ddd$ denotes a sequence of complex numbers of modulus less than one. Individual terms are the values $r_j$ $(j\geq 0)$. 
Continued fractions, whether numerical or functional, have to do with infinite products of $2\times 2$ matrices and infinite composition of the corresponding linear fractional transformations.  (See \cite[\S2.5]{Si:2011} for an exposition of continued fractions from this point of view.)  This paper uses matrix and compositional operator notation for continued fractions, rather than the traditional expressions involving a succession of reciprocals cascading down the page. 

A starting point for the present paper is Verblunsky's version of Szeg\H{o}'s theorem, sketched as follows. 
Work by Schur \cite{Sc:1917,Sc:1918}, Verblunsky \cite{Ve:1936}, Szeg\H{o} \cite{Sz:1975}, Geronimus \cite{Ge:1944,Ge:1946} and others in the first half of the 20th century established a triple correspondence among (i) elements of $\mathcal{H}_1$, (ii) sequences $r:\integer_+\rightarrow\overline{\ddd}$, called Schur parameters, and (iii) probability measures on $\torus$.  (See \cite[Ch.~8]{Kh:2008} for a detailed exposition.)  In particular, if $f\in\mathcal{H}_1$ is neither a finite Blaschke product nor constant, and $f(0)=0$, then the continued fraction expansion of $f$ determines a sequence $r:\integer_+\rightarrow\ddd$ with $r_0=0$, and $f$ determines a probability measure $\mu_f$ on $\torus$ of the form 
\[
d\mu_f=w\,d\tau+d\sigma, 
\]
where $w=\Re(1+f^\circ)/(1-f^\circ)\geq 0$ and $\sigma$ is singular with respect to normalized Lebesgue measure  $\tau$ on the circle. Verblunsky's version of Szeg\H{o}'s theorem asserts 
\begin{equation}\label{verblunsky}
\int_{\torus}\log w\,d\tau=\log\prod_{j=1}^\infty(1-|r_j|^2), 
\end{equation}
with the corollary that 
\begin{equation}\label{gem}
\int_{\torus}\log w\,d\tau>-\infty \Leftrightarrow r\in\ell^2(\integer_+). 
\end{equation}
The latter equivalence has been hailed as a ``gem" of spectral theory \cite[p.~29]{Si:2011}, relating integrability of the spectral quantity $\log w$ to decay of the spatial data $r$.  (See \S\ref{sec-trace}, below, for a physical setting in which $r$ is spatially distributed.)  

For any $d>1$, and any $h\in\mathcal{H}_d$, there is an associated sequence $r:\integer_+\rightarrow\overline{\ddd}$ of Schur parameters determined by the continued fraction expansion of $\diag h$.  Does the analogue of the  Szeg\H{o}-Verblunsky theorem (\ref{verblunsky}) hold in higher dimensions?  Not in general. 
Consider, for example, 
\begin{equation}\label{2d-example}
h(z_1,z_2)=(z_1+z_2)/4\in\mathcal{H}_2.
\end{equation}
The continued fraction expansion of $\diag h(z)=z/2$ yields the sequence $r$, where $r_j=0$ if $j\neq 2$ and $r_2=1/2$ (see \S\ref{sec-continued-fraction} below).  But
\[
\int_{\torus^2}\log\Re\frac{1+h}{1-h}\,d\tau=-\log(112-64\sqrt{3})\neq\log(3/4)=\log\prod_{j=1}^\infty(1-|r_j|^2).
\]
Thus the question becomes
whether for $d>1$ there is a nontrivial subset $\mathcal{S}_d\subset\mathcal{H}_d$ for which the higher dimensional analogue of (\ref{verblunsky}) \emph{does} hold. 
A major objective of the present paper is to characterize explicitly such classes $\mathcal{S}_d$ for every $d\geq1$.

Helson and Lowdenslager derive higher dimensional Szeg\H{o} theorems in \cite[\S2]{HeLo:1958}, although of a different flavour than (\ref{verblunsky})---the sequence of Schur parameters is not involved. Indeed the latter work centres on multivariate Fourier series rather than an approach based on holomorphy. 

Much of the history of Szeg\H{o}'s theorem and its later development, from before Schur's work in 1917 up to the first decade of the present century, is detailed in \cite{Si:2011}.  A second major objective of the present paper is to present a new development that has emerged in the decade since the publication of \cite{Si:2011}, namely the Riemannian structure of the automorphism group of the Poincar\'e disk and its relevance to Schur's original problem of analyzing the Taylor coefficients of functions $f\in\mathcal{H}_1$ \cite{EmSa:2012,BaGi:Proc2017,Gi:JFAA2017}.  Whereas Schur in \cite[p.~210]{Sc:1917} derived a nonlinear recurrence for the Taylor coefficients, the present paper solves this recurrence explicitly, using eigenfunctions of the Laplace-Beltrami operator for the aforementioned Riemannian structure. Precisely the same eigenfunctions represent Taylor coefficients of functions in $\mathcal{H}_d$ that satisfy higher dimensional analogues of the Szeg\H{o}-Verblunsky theorem (\ref{verblunsky}). 

Background on orthogonal polynomials on the unit circle and Szeg\H{o}'s theorem can be found in \cite[Ch.~XI]{Sz:1975}, \cite{Si:1OPUC2005,Si:2005}, \cite[Ch.~8]{Kh:2008} and \cite{Si:2011}.  Many arguments in the present paper rest on the theory of holomorphic and polyharmonic functions on the polydisk; most of the needed results are compiled in Rudin's book \cite{Ru:1969}. See \cite{Ho:1988} for the univariate theory. Some additional background (such as the multivariate version of Montel's theorem) appears in \cite{Oh:2002}.

\subsection{Overview of the paper}

Section~\ref{sec-scattering-disk} introduces the key Riemannian structure on the unit disk, along with the eigenfunctions of its Laplace-Beltrami operator, called scattering polynomials. Scattering weights are defined in \S\ref{sec-weights} as particular tensor products of  scattering polynomials.  

Section~\ref{sec-Schur} defines a closed subset $\mathcal{S}_d\subset\mathcal{H}_d$ of multivariate Schur functions, for each $d\geq 1$. Whereas $\mathcal{S}_1=\mathcal{H}_1$, it turns out that for $d>1$, $\mathcal{S}_d$ is a proper subset of $\mathcal{H}_d$. The continued fraction expansion of Schur functions is analyzed in \S\ref{sec-continued-fraction}. If $d=1$, any given $f\in\mathcal{S}_d$ is determined by its sequence $r:\integer_+\rightarrow\overline{\ddd}$ of Schur parameters generated by its continued fraction expansion.  But for $d>1$ an extra piece of data is required, a sequence \[\nu:\integer_{>0}\rightarrow\{1,\ldots,d\}\] called the variable allocation map.  In general, any Schur function is determined by its corresponding Schur parameters and variable allocation map. Theorem~\ref{thm-full-taylor}, \S\ref{sec-Taylor}, expresses the Taylor coefficients of an arbitrary Schur function in terms of scattering weights evaluated at the Schur parameters.  Section~\ref{sec-standard} introduces the probability measure on the torus determined by a standard Schur function. 

Section~\ref{sec-OPOT} is concerned with orthogonal polynomials on the torus, as determined by a given sequence of Schur parameters and variable allocation map. Section~\ref{sec-construction} presents the basic construction in terms of $2\times2$ matrices; needed properties are derived in \S\ref{sec-properties} and \S\ref{sec-further}; and Theorem~\ref{thm-orthogonality}, \S\ref{sec-orthogonality}, asserts orthogonality of the polynomials constructed in \S\ref{sec-construction} with respect to the corresponding probability measure on the torus. 

Section~\ref{sec-verblunsky} presents three new theorems of Szeg\H{o}-Veblunsky type.  Some technical preliminaries are dispensed with in \S\ref{sec-preliminaries}.  Theorem~\ref{thm-boyd}, \S\ref{sec-multivariate}, is the version needed for later applications; its companion result, Theorem~\ref{thm-multivariate-szego}, is the direct higher dimensional analogue of the classical Szeg\H{o}-Veblunsky theorem (\ref{verblunsky}).  The almost periodic Szeg\H{o}-Veblunsky theorem in \S\ref{sec-almost-periodic}, Theorem~\ref{thm-almost-periodic-main}, applies to the restriction of certain Schur functions to a torus line. This requires a number theoretic result, Proposition~\ref{prop-positive-integer}, as well as Lemma~\ref{lem-continuous-birkhoff} based on Birkhoff's ergodic theorem.  

The underlying motivation for the theorems proved in \S\ref{sec-verblunsky} is an application to one-dimensional scattering theory. This application is presented in the paper's final section~\ref{sec-trace}, in the guise of Theorem~\ref{thm-singular-trace}, a new, singular trace formula for the Schr\"odinger equation in impedance form.  The trace formula is valid for piecewise constant (i.e., discontinuous) impedance functions, a case where the classical trace formula breaks down.

\section{The scattering disk and its orthogonal polynomials\label{sec-scattering-disk}}

\subsection{Irrotational disk automorphisms\label{sec-irrotational}}

Denote by $\mathbb{P}$ the Poincar\'e model of the hyperbolic plane, consisting of $\ddd$ endowed with the metric
\begin{equation}\label{hyperbolic}
ds^2=\frac{4}{(1-x^2-y^2)^2}\left(dx^2+dy^2\right).
\end{equation}
Let $\aut\mathbb{P}$ denote the automorphism group of orientation preserving isometries. Thus $\aut\mathbb{P}\subset\mathcal{H}_1$ consists of all holomorphic bijections of $\ddd$, functions of the form 
\begin{equation}\label{disk-automorphism}
g_{\mu,\rho}:\ddd\rightarrow\ddd,\qquad g_{\mu,\rho}(\xi)=\mu\frac{\xi+\rho}{1+\overline{\rho}\xi}\qquad (\mu\in\torus, \rho\in\ddd, \xi\in\ddd). 
\end{equation}
Call the disk automorphism $g_{\mu,\rho}$ a \emph{rotation} if $\rho=0$, and \emph{irrotational} if $\mu=1$ (equivalently if $g_{\mu,\rho}^\prime(0)>0$).  
Each element of $\aut\mathbb{P}$ factors uniquely as a product of a rotation and an irrotational component,
\begin{equation}\label{factorization}
g_{\mu,\rho}=g_{\mu,0}\circ g_{1,\rho}\quad\mbox{ where }\quad \mu=\frac{g_{\mu,\rho}^\prime(0)}{1-|g_{\mu,\rho}(0)|^2}\quad\mbox{ and }\quad \rho=g_{\mu,\rho}(0)/\mu. 
\end{equation}
The set $\mathfrak{R}$ of rotations is a subgroup of $\aut\mathbb{P}$, but not a normal subgroup.  The set $\mathfrak{I}$ of irrotational elements is not a group; rather, $\mathfrak{I}$ generates $\aut\mathbb{P}$. 
Note that $g_{\mu,\rho}$ extends by the formula (\ref{disk-automorphism}) to a holomorphic function 
\begin{equation}\label{extension}
\tilde{g}_{\mu,\rho}:\overline{\ddd}\rightarrow\overline{\ddd}. 
\end{equation}
(I.e., $\tilde{g}_{\mu,\rho}$ is holomorphic in an open set containing $\overline{\ddd}$.)
Let $g^\circ_{\mu,\rho}$ denote the restriction of this extension to $\torus$. 
By holomorphy, $g^\circ_{\mu,\rho}$  determines $g_{\mu,\rho}$ and vice versa. 
Denote by $\mathbf{c}_r\in\mathcal{H}_1$ the constant function with constant value $r\in\overline{\ddd}$.  The formula (\ref{disk-automorphism}) for $g_{\mu,\rho}$ collapses to the constant function $\mathbf{c}_{\mu\rho}$ if $\rho\in\torus$. 

Remark. Irrotational disk isometries have zero mean rotation $\partial_x\Im g_{1,r}-\partial_y\Re g_{1,r}$ over $\torus$.  More generally, $g_{\mu,r}$ has zero mean rotation over $\mathbb{T}$ only if $\mu=\pm 1$. 

Now, $\aut\mathbb{P}$ carries the topology of uniform convergence on compact sets induced by $\mathcal{H}_1$. With respect to this topology, $\aut\mathbb{P}$ is homeomorphic to the open solid torus via the mapping 
\begin{equation}\label{solid-torus}
\torus\times\ddd\rightarrow\aut\mathbb{P}=\mathfrak{R}\mathfrak{I},\qquad (\mu,\rho)\mapsto g_{\mu,\rho}.
\end{equation}
The closure $\mathfrak{K}$ of $\aut\mathbb{P}$ in $\mathcal{H}_1$ consists of $\aut\mathbb{P}$ together with the circle of constant functions $\mathbf{c}_r$ $(r\in\torus)$. Thus $\mathfrak{K}$ is a semigroup, not a group. Via the map $(\mu,\rho)\mapsto g_{\mu,\rho}$, $\mathfrak{K}$ corresponds to the quotient structure on $\torus\times\overline{\ddd}$ modulo the equivalence relation 
\begin{equation}\label{equivalence-relation}
(\mu_1,\rho_1)\sim(\mu_2,\rho_2)\quad\Longleftrightarrow\quad (\mu_1,\rho_1)=(\mu_2,\rho_2)\quad\mbox{ or }\quad\mu_1\rho_1=\mu_2\rho_2\in\torus. 
\end{equation}
The resulting structure, which in effect adjoins a circle to $\torus\times\ddd$, is homeomorphic to $S^3$,
\begin{equation}\label{S3}
\torus\times\overline{\ddd}/\sim\;\;\cong\mathfrak{K}=\aut\mathbb{P}\cup\left\{\mathbf{c}_r\,|\,r\in\torus\right\}\cong S^3.
\end{equation}

The important part of this picture for present purposes is the closure in $\mathfrak{K}$ of the set $\mathfrak{I}$ of irrotational disk automorphisms, corresponding to the closed disk $\overline{\ddd}$.  This closed disk carries a natural Riemannian structure (degenerate at the boundary circle) that yields an explicit solution to Schur's original 1917 problem of characterizing the Taylor coefficients associated to elements of $\mathcal{H}_1$.

\subsection{The scattering disk\label{sec-scattering}}

Emmanuele and Salvai in \cite{EmSa:2012} derived a Riemannian structure on $\aut\mathbb{P}$ by considering the action of $g^\circ_{\mu,\rho}$ on $\torus$, and minimizing the total kinetic energy on $\torus$ along geodesics in $\aut\mathbb{P}$. The result is a product metric corresponding to the factorization $\aut\mathbb{P}=\mathfrak{R}\mathfrak{I}\cong\torus\times\ddd$. The factor on the rotational component $\mathfrak{R}\cong\torus$ is the standard Euclidean metric on the circle; the factor on the irrotational component $\mathfrak{I}\cong\ddd$ is what we call the \emph{scattering metric},
\begin{equation}\label{scattering-metric}
ds^2=\frac{4}{1-x^2-y^2}\left(dx^2+dy^2\right). 
\end{equation}
The term scattering metric stems from independent work \cite{Gi:JFAA2017,Gi:JAT2019}, where the same metric arose in the context of waves in layered media. The latter work derives the eigenfunctions of the Laplace-Beltrami operator for (\ref{scattering-metric}), 
\begin{equation}\label{laplace-beltrami}
\Delta=-\frac{1-x^2-y^2}{4}\left(\frac{\partial^2}{\partial x^2}+\frac{\partial^2}{\partial y^2}\right)=-(1-z\bar{z})\frac{\partial^2}{\partial\bar{z}\partial z}\qquad (z=x+iy).
\end{equation}
Remarkably, the eigenfunctions turn out to be bivariate polynomials, orthogonal with respect to the area measure 
\[
\frac{4}{1-x^2-y^2}dxdy
\]
determined by (\ref{scattering-metric}). These \emph{scattering polynomials} are central to present considerations.  

Much more is true about the scattering disk. Indeed, the unit disk endowed with the metric (\ref{scattering-metric}) has so many miraculous properties it seems hard to believe it was not known to the pioneers of differential geometry; yet apparently these properties were discovered only recently. 
Here is one illustrative example. The eigenvalues of (\ref{laplace-beltrami}) are the positive integers, and the multiplicity of $n$ as an eigenvalue is the number of divisors of $n$. It follows that the spectral zeta function for the scattering disk is the square of the Riemann zeta function, and so has the same zeros. See \cite{BaGi:Proc2017}.

\subsection{Eigenfunctions of the Laplace-Beltrami operator\label{sec-scattering-polynomials}}
The definition and key properties are as follows. 
\begin{defn}\label{defn-scattering-polynomials}
For each $(p,q)\in\integer_+^2$ define $\varphi^{(p,q)}:\complex\rightarrow\complex$ as follows. If $\min\{p,q\}\geq 1$ set
\begin{equation}\label{varphi}
\varphi^{(p,q)}(z)=\textstyle\frac{(-1)^p}{q(p+q-1)!}\,\displaystyle(1-z\bar{z})\frac{\partial^{\,p+q}}{\partial\bar{z}^p\partial z^q}(1-z\bar{z})^{p+q-1}.
\end{equation}
If $p=0<q$ set $\varphi^{(p,q)}=0$, and if $p\geq 0$ set $\varphi^{(p,0)}(z)=z^p$.  The functions $\varphi^{(p,q)}$ so defined are polynomials with respect to variables $z$ and $\bar{z}$ (or $x$ and $y$, where $z=x+iy$), referred to as \emph{scattering polynomials}. 
\end{defn}
In evaluating $\varphi^{(p,q)}(0)$, the intended interpretation is that $z^0=1$ even if $z=0$. For example, $\varphi^{(0,0)}(z)=1$ is constant.
\begin{prop}[{From \cite[Thm.~1.3]{BaGi:Proc2017} and \cite[p.~39]{Gi:JAT2019}}]\label{prop-scattering-props}
For every $(p,q)\in\integer_+^2$,
\begin{equation}\label{eigenvalue}
\Delta\varphi^{(p,q)}=pq\varphi^{(p,q)},
\end{equation}
where $\Delta$ is given by (\ref{laplace-beltrami}), and every eigenfunction of $\Delta$ is a scattering polynomial. If $\min\{p,q\}\geq 1$, then
\begin{equation}\label{explicit}
\varphi^{(p,q)}(z)=\sum_{j=1}^{\min\{p,q\}}\binom{p}{j}\binom{q-1}{j-1}z^{p-j}(-\bar{z})^{q-j}(1-z\bar{z})^j.
\end{equation}
\end{prop}
Remarks. The proposition makes obvious the multiplicity of eigenvalues, and shows that scattering polynomials have integer coefficients.

\subsection{Scattering weights\label{sec-weights}}

\begin{defn}\label{defn-scattering-weight}
Let $\alpha:\integer_{>0}\rightarrow\integer_+$ have finite support (i.e. be eventually zero). Set $n=0$ if $\alpha$ is identically zero; otherwise set $n=\max\supp\alpha$.  Given $r:\integer_+\rightarrow\ddd$, set 
\begin{equation}\label{scattering-weight}
c_\alpha(r)=\varphi^{(1,\alpha_1)}(r_0)\prod_{j=1}^n\varphi^{(\alpha_j,\alpha_{j+1})}(r_j).
\end{equation}
Call $c_\alpha$ the \emph{scattering weight} associated with multi-index $\alpha$. 
\end{defn}
Remarks. By Definition~\ref{defn-scattering-polynomials} and Proposition~\ref{prop-scattering-props}, 
\begin{equation}\label{zero-values}
\varphi^{(\alpha_j,\alpha_{j+1})}(r_j)=0\quad\mbox{ if }\quad\alpha_j=0\neq\alpha_{j+1}, 
\end{equation}
irrespective of $r_j$; and 
\begin{equation}\label{phi-r0}
\varphi^{(1,\alpha_1)}(0)=\left\{\begin{array}{cc}1&\mbox{ if }\alpha_1=1\\
0&\mbox{ otherwise }
\end{array}\right..
\end{equation}
Thus if $r_0=0$, then $c_\alpha$ is non-zero only if $\alpha_1=1$. 
If $n=\max\supp\alpha>0$, then (\ref{zero-values}) implies $c_\alpha$ is non-zero only if $\alpha$ has contiguous support, i.e., only if $\supp\alpha=\{1,\ldots,n\}$. 

\section{Multivariate Schur functions\label{sec-Schur}}

Given $1\leq j\leq d$, define
\begin{equation}\label{multiplication}
\mathfrak{m}_j:\mathcal{H}_d\rightarrow\mathcal{H}_d,\qquad (\mathfrak{m}_jh)(z)=z_jh(z)\qquad\bigl(h\in\mathcal{H}_d,\; z\in\ddd^d\bigr).
\end{equation}
Recalling the notation $\tilde{g}_{\mu,r}$ from (\ref{extension}), for any $r\in\overline{\ddd}$ and $d\geq 1$, define
\begin{equation}\label{composition}
\mathfrak{g}_r:\mathcal{H}_d\rightarrow\mathcal{H}_d,\qquad\mathfrak{g}_rh=\tilde{g}_{1,r}\circ h\qquad\bigl(h\in\mathcal{H}_d\bigr).
\end{equation}
Remark. It should be emphasized that $\mathfrak{g}_r$ is a \emph{nonlinear} operator on $\mathcal{H}_d$, except in the trivial case $r=0$, whereas the multiplication operator $\mathfrak{m}_j$ is linear. 

If $r\in\torus$ then the operator $\mathfrak{g}_r$ is continuous on $\mathcal{H}_d$, since it is constant; i.e., $\mathfrak{g}_rh=\mathbf{c}_r$ for every $h\in\mathcal{H}_d$. 
For any $1\leq j\leq d$ and $r\in\ddd$, $\mathfrak{m}_j$ and $\mathfrak{g}_r$ are homeomorphisms on $\mathcal{H}_d$, since for any $h_1,h_2\in\mathcal{H}_d$, 
\begin{equation}\label{m-homeomorphism}
\left|\mathfrak{m}_jh_1-\mathfrak{m}_jh_2\right|=\left|h_1-h_2\right|
\end{equation}
and 
\begin{equation}\label{r-homeomorphism}
\frac{1-|r|}{1+|r|}\left|\mathfrak{g}_rh_1-\mathfrak{g}_rh_2\right|\leq\left|h_1-h_2\right|
\leq \frac{1+|r|}{1-|r|}\left|\mathfrak{g}_rh_1-\mathfrak{g}_rh_2\right|.
\end{equation}
Note also that for any $r_1,r_2\in\ddd$ and $h\in\mathcal{H}_d$,
\begin{equation}\label{gr-weak-continuity}
\left|\mathfrak{g}_{r_1}h(z)-\mathfrak{g}_{r_2}h(z)\right|\leq\frac{4|r_1-r_2|}{(1-|r_1|)(1-|r_2|)}\qquad\bigl(z\in\ddd^d\bigr).
\end{equation}
The operator $\mathfrak{g}_r$ is therefore \emph{weakly continuous} with respect to $r\in\ddd$, in the sense that if $r_n\rightarrow r$, where $r,r_n\in\ddd$, then for any $h\in\mathcal{H}_d$, 
\[
\mathfrak{g}_{r_n}h\rightarrow\mathfrak{g}_rh
\]
uniformly on $\ddd^d$, by (\ref{gr-weak-continuity}).  And if $r\in\torus$, in which case $\mathfrak{g}_rh=\mathbf{c}_r$, one can show $\mathfrak{g}_{r_n}h\rightarrow\mathbf{c}_r$ uniformly on compact subsets of $\ddd^d$, provided $h$ is not itself the constant function $\mathbf{c}_{-r}$. (The case $r_n\rightarrow r\in\torus$ is not needed in the sequel.)

For each $d\geq1$, let $\mathfrak{S}_d$ denote the semigroup generated by (nonlinear) product-composition operators on $\mathcal{H}_d$ of the form $\mathfrak{g}_r\mathfrak{m}_j$ where $1\leq j\leq d$ and $r\in\ddd$. Let $\mathfrak{S}_d\mathbf{c}_0\subset\mathcal{H}_d$ denote the orbit of the zero function under the action of $\mathfrak{S}_d$. 
\begin{defn}\label{defn-schur}
Define $\mathcal{S}_d$ to be the closure in $\mathcal{H}_d$ of the set $\mathfrak{S}_d\mathbf{c}_0$.  Elements of $\mathcal{S}_d$ will be referred to as \emph{Schur functions}.  
\end{defn}

Remarks. 
Schur himself considered only the case $d=1$, proving that $\mathcal{S}_1=\mathcal{H}_1$. In higher dimensions $\mathcal{S}_d$ is a proper  subset of $\mathcal{H}_d$ (see (\ref{2d-example}), above). Schur functions are by definition limits of sequences of terms of the form
\begin{equation}\label{term}
\mathfrak{g}_{r_0}\mathfrak{m}_{\nu_1}\mathfrak{g}_{r_1}\mathfrak{m}_{\nu_2}\cdots\mathfrak{g}_{r_{k-1}}\mathfrak{m}_{\nu_k}\mathbf{c}_0.
\end{equation}
Note that
since $\mathbf{c}_0=\mathfrak{g}_0\mathfrak{m}_1\mathbf{c}_0$, one is free to choose a representation (\ref{term}) in which the terminal index $k$ is as large as desired.  Moreover, the form (\ref{term}) is equivalent to 
\begin{equation}\label{term-2}
\mathfrak{g}_{r_0}\mathfrak{m}_{\nu_1}\mathfrak{g}_{r_1}\mathfrak{m}_{\nu_2}\cdots\mathfrak{g}_{r_{k-1}}\mathbf{c}_0,
\end{equation}
since $\mathfrak{m}_{\nu_k}\mathbf{c}_0=\mathbf{c}_0$.

\subsection{Continued fractions of Schur functions\label{sec-continued-fraction}}

Every Schur function $f\in\mathcal{S}_d$ can be unfolded as a continued fraction; the process determines two sequences, a sequence of Schur parameters and a variable allocation map, which in turn uniquely characterize $f$.   

\begin{lem}\label{lem-continued-fraction}
Let $f\in\mathcal{S}_d$ be arbitrary, and set $r=f(0)$.  There exist $1\leq j\leq d$ and $h\in\mathcal{S}_d$ such that 
\[
f=\mathfrak{g}_r\mathfrak{m}_jh. 
\]
\end{lem}
\begin{pf} If $r\in\torus$, then $f=\mathbf{c}_r$ by the maximum principle. In this case,
\[
f=\mathfrak{g}_r\mathfrak{m}_1\mathbf{c}_0,
\]
where $\mathbf{c}_0\in\mathcal{S}_d$. 
Suppose on the other hand that $r\in\ddd$.  
Since $f\in\mathcal{S}_d$, there exists a sequence $f_n\in\mathfrak{S}_d\mathbf{c}_0$ such that $f_n\rightarrow f$ in $\mathcal{H}_d$. Write
\[
f_n=\mathfrak{g}_{r_{0,n}}\mathfrak{m}_{\nu_{1,n}}\cdots\mathfrak{g}_{r_{k_n-1,n}}\mathfrak{m}_{\nu_{k_n,n}}\mathbf{c}_0\in\mathfrak{S}_d\mathbf{c}_0\qquad (n\geq1),
\]
where, without loss of generality, each $k_n\geq 2$. Note that the sequence of indices $\nu_{1,n}$ $(n\geq 1)$ takes at least one value $j\in\{1,\ldots,d\}$ infinitely often; passing to the corresponding subsequence $f_{m_n}$ $(n\geq 1)$ yields a sequence for which $\nu_{1,m_n}=j$ is constant. Thus no loss of generality results from assuming $\nu_{1,n}=j$ $(n\geq 1)$ is constant to begin with, whereby
\[
f_n=\mathfrak{g}_{r_{0,n}}\mathfrak{m}_jh_n\quad\mbox{ with }\quad h_n=
\mathfrak{g}_{r_{1,n}}\mathfrak{m}_{\nu_{2,n}}\cdots\mathfrak{g}_{r_{k_n-1,n}}\mathfrak{m}_{\nu_{k_n,n}}\mathbf{c}_0
\in\mathfrak{S}_d\mathbf{c}_0\qquad(n\geq1).
\]
Observe that
\begin{equation}\label{r-approximation}
\begin{split}
|\mathfrak{g}_r\mathfrak{m}_jh_n-f|&\leq |\mathfrak{g}_r\mathfrak{m}_jh_n-f_n|+|f_n-f|\\
&=|\mathfrak{g}_r\mathfrak{m}_jh_n-\mathfrak{g}_{r_{0,n}}\mathfrak{m}_jh_n|+|f_n-f|\\
&\leq\frac{4|r-r_{0,n}|}{(1-|r|)(1-|r_{0,n}|)}+|f_n-f|\quad\mbox{ by (\ref{gr-weak-continuity})}.
\end{split}
\end{equation}
Now,
\[
r_{0,n}=f_n(0)\rightarrow f(0)=r,
\]
since $f_n\rightarrow f$.  
Therefore (\ref{r-approximation}) implies
\begin{equation}\label{intermediate-convergence-step-1}
\mathfrak{g}_r\mathfrak{m}_jh_n\rightarrow f
\end{equation}
uniformly on compact subsets of $\ddd^d$.  
For every $m,n\geq 1$, 
\[
\begin{split}
|h_m-h_n|&=|\mathfrak{m}_jh_m-\mathfrak{m}_jh_n| \\
&\leq \frac{1+|r|}{1-|r|}\left|\mathfrak{g}_r\mathfrak{m}_jh_m-\mathfrak{g}_r\mathfrak{m}_jh_n\right|\quad\mbox{ by (\ref{r-homeomorphism})}.
\end{split}
\]
Thus by (\ref{intermediate-convergence-step-1}) the sequence $h_n\in\mathfrak{S}_d\mathbf{c}_0$ is uniformly Cauchy on compact subsets of $\ddd^d$, and hence convergent in $\mathcal{H}_d$ to some $h\in\mathcal{S}_d$. It follows in turn by continuity of $\mathfrak{g}_r\mathfrak{m}_j$ that 
\[
\mathfrak{g}_r\mathfrak{m}_jh_n\rightarrow\mathfrak{g}_r\mathfrak{m}_jh=f.
\]
\end{pf}

Remark.  The index $j$ is not necessarily uniquely determined, since
whenever $r_j=0$, commutativity of the two adjacent multiplication operators allows the order to be reversed:
\[
\mathfrak{m}_j\mathfrak{g}_0\mathfrak{m}_k=\mathfrak{m}_j\mathfrak{m}_k=\mathfrak{m}_k\mathfrak{g}_0\mathfrak{m}_j.  
\]
(The non-uniqueness could be eliminated by defining a reduced representation in the obvious way.) 
Lemma~\ref{lem-continued-fraction} justifies the following construction. 
\begin{defn}\label{defn-continued-fraction}
Let $f\in\mathcal{S}_d$. Construct sequences 
$r_n\in\overline{\ddd}$ $(n\geq 0)$, $\nu_n\in\{1,\ldots,d\}$ $(n\geq 1)$ and 
$h_n\in\mathcal{S}_d$ $(n\geq 1)$ by starting with $h_0=f$ and iterating the following. 
Given $h_n\in\mathcal{S}_d$ for some $n\geq 0$, set $r_n=h_n(0)$, and apply
Lemma~\ref{lem-continued-fraction} to obtain $\nu_{n+1}\in\{1,\ldots,d\}$ and $h_{n+1}\in\mathcal{S}_d$ such that 
\begin{equation}\label{iteration}
h_n=\mathfrak{g}_{r_n}\mathfrak{m}_{\nu_{n+1}}h_{n+1}.  
\end{equation}
In particular, if $h_n$ is constant, set $\nu_{n+1}=1$ and $h_{n+1}=\mathbf{c}_0$. The corresponding terms 
\begin{equation}\label{convergents}
f_n:=\mathfrak{g}_{r_0}\mathfrak{m}_{\nu_1}\mathfrak{g}_{r_1}\cdots\mathfrak{m}_{\nu_n}\mathfrak{g}_{r_n}\mathbf{c}_0\quad\mbox{ and }\quad k_n:=\mathfrak{m}_{\nu_{n+1}}h_{n+1}\qquad(n\geq0)
\end{equation} 
are respectively called \emph{convergents} and \emph{fractional parts} in the continued fraction expansion of $f$.  The numbers $r_n$ $(n\geq 0)$ are called the \emph{Schur parameters} of $f$, and 
\[
\nu:\integer_{>0}\rightarrow\{1,\ldots,d\}
\]
is called the \emph{variable allocation} map.  Denote by $\mathfrak{f}_n$ the operator associated with the $n$th convergent,
\begin{equation}\label{frak-fn}
\mathfrak{f}_n:=\mathfrak{g}_{r_0}\mathfrak{m}_{\nu_1}\mathfrak{g}_{r_1}\cdots\mathfrak{m}_{\nu_n}\mathfrak{g}_{r_n}\qquad(n\geq 0),
\end{equation}
so that $f$ has the representations
\begin{equation}\label{f-factorizations}
f=\mathfrak{f}_nk_n\qquad(n\geq 0). 
\end{equation}
\end{defn}
Remarks. The sequence $r$ arising in Definition~\ref{defn-continued-fraction} is uniquely determined by $f$, while the variable allocation map $\nu$ is not.  Observe however that different possible choices of $\nu$ all give rise to the same sequence of convergents; i.e. the functions $f_n$ are uniquely determined.  Note also that $f_n,k_n\in\mathcal{S}_d$, and $k_n(0)=0$.

\subsection{Taylor series of Schur functions\label{sec-Taylor}}

Let $A_0=\{\mathbf{0}\}$, where 
\[
\mathbf{0}:\integer_{>0}\rightarrow\integer_+,\qquad\mathbf{0}_n=0\qquad(n\geq 1)
\]
denotes the zero multi-index. For each $n\geq 1$, let $A_n$ denote the set of all multi-indices $\alpha:\integer_{>0}\rightarrow\integer_+$ such that 
$
\supp\alpha\subset\{1,\ldots,n\}, 
$
and set 
\begin{equation}\label{A}
A=\bigcup_{n=0}^\infty A_n.
\end{equation}
It follows from (\ref{explicit}) that if $r_n\in\torus$ then $c_\alpha(r)=0$ unless $\alpha\in A_n$, for every $n\geq0$.  Note also that if $\alpha\in A_n$ and $r:\integer_+\rightarrow\overline{\ddd}$,  then $c_\alpha(r)$ depends only on $r_j$ for $j\leq n$. 

Given $\nu:\integer_{>0}\rightarrow\{1,\ldots,d\}$ and $\alpha\in A_n$, let $(z\circ\nu)^\alpha$ denote the monomial
\begin{equation}\label{composition-notation}
(z\circ\nu)^\alpha=z_{\nu_1}^{\alpha_1}\cdots z_{\nu_n}^{\alpha_n}.
\end{equation}

\begin{prop}\label{prop-convergent-representation}
Let $c_\alpha$ denote the scattering weights of Definition~\ref{defn-scattering-weight}.
Let $f_n\in\mathfrak{S}_d\mathbf{c}_0$ have the form
\begin{equation}\label{fn-form}
f_n=\mathfrak{g}_{r_0}\mathfrak{m}_{\nu_1}\cdots\mathfrak{g}_{r_n}\mathfrak{m}_{\nu_{n+1}}\mathbf{c}_0,
\end{equation}
where $n\geq 0$, and let $r:\integer_+\rightarrow\overline{\ddd}$ and $\nu:\integer_{>0}\rightarrow\{1,\ldots,d\}$ be arbitrary extensions of the initial sequences $r_0,\ldots,r_n$ and $\nu_1,\ldots,\nu_{n+1}$ occurring in (\ref{fn-form}). 
Then the Taylor series expansion of $f_n$ is
\begin{equation}\label{finite-support-taylor}
f_n(z)=\sum_{\alpha\in A_n}c_\alpha(r)(z\circ\nu)^\alpha,
\end{equation}
which converges absolutely and uniformly on $\overline\ddd^d$. 
\end{prop}
\begin{pf}
To begin, set 
\[
g_n=\mathfrak{g}_{r_0}\mathfrak{m}_{1}\cdots\mathfrak{g}_{r_n}\mathfrak{m}_{{n+1}}\mathbf{c}_0,
\]
so that $f_n(z)=g_n(z\circ\nu)$. 
Introduce an auxiliary variable $z_0$ and corresponding multiplication operator $\mathfrak{m}_0$. Set 
\[
g^{\dagger}_n=\mathfrak{m}_0g_n. 
\]
Note that $g^{\dagger}_n$ is a holomorphic rational function with respect to $(z_0,\ldots,z_n)\in\overline{\ddd}^{n+1}$.  It is proved in \cite[Thm.~1]{Gi:JFAA2017} that the restriction of $g^{\dagger}_n$ to $\torus^{n+1}$ is represented by the pointwise convergent Fourier series
\begin{equation}\label{fourier-series-0}
g^{\dagger}_n(z)=\sum_{\alpha\in\{1\}\times A_n}\left(\prod_{j=0}^n\varphi^{(\alpha_j,\alpha_{j+1})}(r_j)\right)z^\alpha\qquad(z\in\torus^{n+1}).
\end{equation}
One recovers $g_n$ from $g^{\dagger}_n$ by restricting to $z_0=1$, whereupon the formula (\ref{fourier-series-0}) reduces to 
\begin{equation}\label{fourier-series-reduced}
g_n(z)=\sum_{\alpha\in A_n}c_\alpha(r)z^\alpha\qquad(z\in\torus^n).
\end{equation}
Since the domain of holomorphy of $g_n$ includes an open polydisk containing $\overline{\ddd}^n$, and since such an open polydisk is a complete Reinhardt domain, it follows that the series (\ref{fourier-series-reduced}) converges absolutely and uniformly on $\overline{\ddd}^n$.  But then the corresponding series for $f_n(z)=g_n(z\circ\nu)$ converges absolutely and uniformly on $\overline{\ddd}^d$. 
\end{pf}

\begin{prop}\label{prop-h-form}
Let $f\in\mathcal{S}_d$ have the form 
\begin{equation}\label{f-h-form}
f=\mathfrak{g}_{r_0}\mathfrak{m}_{\nu_1}\cdots\mathfrak{g}_{r_n}\mathfrak{m}_{\nu_{n+1}}h,
\end{equation}
where $n\geq 0$ and $h\in\mathcal{S}_d$. Let $\nu:\integer_{>0}\rightarrow\{1,\ldots,d\}$ be an arbitrary extension of the initial sequence $\nu_1,\ldots,\nu_{n+1}$ occurring in (\ref{f-h-form}). 
Then $f$ has the representation
\begin{equation}\label{finite-support-taylor-2}
f(z)=\sum_{\alpha\in A_{n+1}}\varphi^{(1,\alpha_1)}(r_0)\left(\prod_{j=1}^n\varphi^{(\alpha_j,\alpha_{j+1})}(r_j)\right)h(z)^{\alpha_{n+1}}(z\circ\nu)^\alpha,
\end{equation}
which converges absolutely and uniformly on $\ddd^d$. 
\end{prop}
\begin{pf}
Introduce an auxiliary variable $\xi\in\complex$, let $r:\integer_+\rightarrow\overline{\ddd}$ be an arbitrary extension of the initial sequence $r_0,\ldots,r_n,\xi$, and set 
\[
g(z;\xi)=\mathfrak{g}_{r_0}\mathfrak{m}_{\nu_1}\cdots\mathfrak{g}_{r_n}\mathfrak{m}_{\nu_{n+1}}\mathfrak{g}_\xi\mathbf{c}_0.
\]
Then for each fixed $\xi\in\overline{\ddd}$, Proposition~\ref{prop-convergent-representation} implies
\begin{equation}\label{g-z-xi}
g(z;\xi)=\sum_{\alpha\in A_{n+1}}c_\alpha(r)(z\circ\nu)^\alpha,
\end{equation}
with the series converging absolutely and uniformly for $z\in\ddd^d$. 
Note that for $\alpha\in A_{n+1}$ (and with $r_{n+1}=\xi$),
\[
c_\alpha(r)=\varphi^{(1,\alpha_1)}(r_0)\left(\prod_{j=1}^n\varphi^{(\alpha_j,\alpha_{j+1})}(r_j)\right)\xi^{\alpha_{n+1}},
\]
whence
\[
g(z;\xi)=\varphi^{(1,\alpha_1)}(r_0)\left(\prod_{j=1}^n\varphi^{(\alpha_j,\alpha_{j+1})}(r_j)\right)\xi^{\alpha_{n+1}}(z\circ\nu)^\alpha.
\]
Absolute and uniform convergence of the latter series relative to $z\in\overline{\ddd}^d$ for a fixed $\xi\in\torus$ implies absolute and uniform convergence relative to $(z,\xi)\in\overline{\ddd}^{d+1}$.  The desired result then follows from the facts that $f(z)=g(z;h(z))$ and $h:\ddd^d\rightarrow\overline{\ddd}$. 
\end{pf}

For each $n\geq1$, let $\pi_n:\mathcal{H}_d\rightarrow\mathcal{H}_d$ denote projection onto the degree $n$ Taylor polynomial.  It follows from the multivariate Cauchy formula that each $\pi_n$ is a continuous operator.
The foregoing two propositions yield the following.
\begin{lem}\label{lem-convergent-taylor}
Let $f\in\mathcal{S}_d$ have $n$th convergent $f_n$ as per Definition~\ref{defn-continued-fraction}.  Then $\pi_nf=\pi_nf_n$. 
\end{lem}
\begin{pf}
If $|\alpha|\leq n$ and either $c_\alpha(r)\neq 0$ or, referring to (\ref{finite-support-taylor-2}), 
\[
\varphi^{(1,\alpha_1)}(r_0)\left(\prod_{j=1}^n\varphi^{(\alpha_j,\alpha_{j+1})}(r_j)\right)\neq 0,
\]
then $\alpha_{n+1}=0$ and $\alpha\in A_n$.  Therefore the coefficient of $(z\circ\nu)^\alpha$ in (\ref{finite-support-taylor-2}) is 
\[
\varphi^{(1,\alpha_1)}(r_0)\left(\prod_{j=1}^n\varphi^{(\alpha_j,\alpha_{j+1})}(r_j)\right)=c_\alpha(r),
\]
the same as in (\ref{finite-support-taylor}). 
\end{pf}

\begin{thm}\label{thm-full-taylor}
Let $f\in\mathcal{S}_d$ be given. Let $r:\integer_+\rightarrow\overline{\ddd}$ and $\nu:\integer_{>0}\rightarrow\{1,\ldots,d\}$ be the Schur parameters and variable allocation map, as per Definition~\ref{defn-continued-fraction}, and let $f_n$ $(n\geq 1)$ be the sequence of convergents. 
Then 
$f_n\rightarrow f$ uniformly on compact subsets of $\ddd^d$, and 
\begin{equation}\label{full-taylor}
f(z)=\sum_{\alpha\in A}c_\alpha(r)(z\circ\nu)^\alpha\qquad \bigl(z\in\ddd^d\bigr),
\end{equation}
with the latter series converging uniformly on compact sets in $\ddd^d$. 
\end{thm}
\begin{pf}  Observe that if $c_\alpha(r)\neq0$ and $|\alpha|\leq n$, then $\alpha\in A_n$.  It follows by Proposition~\ref{prop-convergent-representation} that for any fixed $m\geq 1$ the sequence $\pi_mf_n$ $(n\geq 1)$ is eventually constant, with constant value 
\begin{equation}\label{taylor-polynomial}
\pi_mf_m(z)=\sum_{\stackrel{\alpha\in A_m}{|\alpha|\leq m}}c_\alpha(r)(z\circ\nu)^\alpha\qquad \bigl(z\in\ddd^d\bigr).
\end{equation}
The multivariate version of Montel's theorem \cite[Thm.~1.5]{Oh:2002} implies that $\mathcal{H}_d$ is compact; in particular, $f_n$ $(n\geq 1)$ has a convergent subsequence $f_{j_n}$ $(n\geq 1)$ converging to some $g\in\mathcal{H}_d$.  The degree $m$ Taylor polynomial of $g$ is given by the right-hand side of (\ref{taylor-polynomial}), and therefore $g$ has Taylor expansion given by the right-hand side of (\ref{full-taylor}).  It follows that $g$ is the unique limit point of the sequence $f_n$ $(n\geq 1)$, and hence that $f_n\rightarrow g$.  To prove that $g=f$ it suffices to note by Lemma~\ref{lem-convergent-taylor} that 
\[
\pi_nf=\pi_nf_n=\pi_ng, 
\]
whence $f$ and $g$ have the same Taylor series. 
\end{pf}

\subsection{Standard Schur functions and their probability distributions\label{sec-standard}}

Recall the notation $\tau$ to denote normalized Lebesgue measure on $\torus^d$, so that, for $f\in L^1(\torus^d)$, 
\begin{equation}\label{tau-notation}
\int_{\torus^d}f\,d\tau=\frac{1}{(2\pi)^d}\displaystyle\rule{-27pt}{0pt}\int\limits_{(e^{i\theta_1},\ldots,e^{i\theta_d})\in\torus^d}\rule{-27pt}{0pt}f(e^{i\theta_1},\ldots,e^{i\theta_d})\,d\theta_1\cdots d\theta_d.
\end{equation} 
The same symbol $\tau$ is used for all dimensions, with the pertinent value of $d$ being evident from context. 
\begin{defn}\label{defn-standard-schur}
A Schur function $f\in\mathcal{S}_d$ is termed \emph{standard} if $f(0)=0$ and each Schur parameter $r_n$ has modulus strictly less than 1 $(n\geq0)$. 
\end{defn}
Remark. A standard Schur function $f$ takes values in $\ddd$. For, a function in $\mathcal{H}_d$ takes a value on $\torus$ only if it is constant, by the maximum principle, whereas $f(0)=0\not\in\torus$. 

A standard Schur function $f$ has an associated probability measure $\mu_f$ on the torus $\torus^d$, defined as follows.
Given $z=(s_1e^{i\phi_1},\ldots,s_de^{i\phi_d})\in\ddd^d$, denote by $K_z\in C(\torus^d)$ the associated Poisson kernel
\begin{equation}\label{poisson-kernel}
K_z(e^{i\theta_1},\ldots,e^{i\theta_d})=\prod_{j=1}^d\frac{1-s_j^2}{1-2s_j\cos(\phi_j-\theta_j)+s_j^2}\qquad\bigl((e^{i\theta_1},\ldots,e^{i\theta_d})\in\torus^d\bigr).
\end{equation}
As in the single-variable case, under certain conditions a $d$-harmonic function on $\ddd^d$ can be recovered from integration of its radial limit against $K_z$.  

\begin{prop}[{from \cite[Ch.~2,3]{Ru:1969}}]\label{prop-poisson-facts}
\begin{enumerate}[label={(\roman*)},itemindent=0pt]
\item\label{bounded-poisson}For any $h\in\mathcal{H}_d$, 
the radial limit 
\[
h^\circ(z)=\lim_{s\nearrow 1}h(sz)
\]
exists for a.e.~$z\in\torus^d$, and 
\[
h(z)=\int_{\torus^d}K_zh^\circ\,d\tau\qquad(z\in\ddd^d).
\]
\item\label{weak-boundary-convergence} If $h_n\rightarrow h$ in $\mathcal{H}_d$, then $h_n^\circ\rightarrow h^\circ$ weakly in the dual space to $C(\torus^d)$.\\[0pt]
\item\label{associated-measure}
For every $f\in\mathcal{H}_d$ with $f(0)=0$, there is a unique measure $\mu_f$ on $\torus^d$ of the form 
\[
d\mu_f=w\,d\tau+d\sigma\quad\mbox{ where }\quad w=\Re\frac{1+f^\circ}{1-f^\circ}
\]
and $\sigma$ is singular with respect to Lebesgue measure on $\torus^d$, such that 
\[
\Re\frac{1+f(z)}{1-f(z)}=\int_{\torus^d}K_z\,d\mu_f\qquad(z\in\ddd^d). 
\]
\item\label{continuous-boundary-function}
If $f\in\mathcal{H}_d$ extends to a holomorphic function $\tilde{f}:\overline{\ddd}^d\rightarrow\ddd$ and $f(0)=0$, then $\mu_f$ is absolutely continuous with respect to $\tau$, and $w=d\mu_f/d\tau\in C(\torus^d)$. 
\end{enumerate}
\end{prop}
\begin{pf}
(i) This is the multivariate version of Fatou's theorem; see \cite[\S2.3.2(a)]{Ru:1969}. 

(ii) By Alaoglu's theorem the sequence $h_n^\circ$ has a weak limit point $G$ in the dual space to $C(\torus^d)$.  If $h_{j_n}\rightarrow G$ weakly as $n\rightarrow\infty$, then for every $z\in\ddd^d$, 
\[
G(K_z)=\lim_{n\rightarrow\infty}\int_{\torus^d}K_zh_{j_n}^\circ\,d\tau=\lim_{n\rightarrow\infty}h_{j_n}(z)=h(z)=
\int_{\torus^d}K_zh^\circ\,d\tau. 
\]
Thus $G$ is uniquely determined and (identifying a function with its induced functional) $G=h^\circ$, by the Riesz representation theorem; moreover, $h_n^\circ\rightarrow h^\circ$ weakly since $h^\circ$ is the unique limit point of the sequence $h_n^\circ$. 

(iii) Consider the positive, $d$-harmonic function
\begin{equation}\label{k-real-part}
k(z)=\Re\frac{1+f(z)}{1-f(z)}\qquad(z\in\ddd^d). 
\end{equation}
By the mean value property of polyharmonic functions, for each $0\leq s<1$, 
\begin{equation}\label{k-mean-values}
\int_{\torus^d}k(sz)\,d\tau(z)=s^dk(0)=s^d<1. 
\end{equation}
Thus the family $k(s\,\cdot)$ is uniformly bounded in $L^1(\torus^d)$, implying the existence of a unique positive Borel measure $\mu_f$ on $\torus^d$ such that for every $z\in\ddd^d$,
\begin{equation}\label{reproducing-measure}\nonumber
k(z)=\int_{\torus^d}K_z\,d\mu_f.
\end{equation}
See \cite[Thm.~2.1.3(e)]{Ru:1969}.  In particular,
\[
1=k(0)=\int_{\torus^d}K_0\,d\mu_f=\int_{\torus^d}d\mu_f,
\]
whence $\mu_f$ is a probability measure.  Moreover, by \cite[Thm.~3.2.4(iv)]{Ru:1969}, the condition (\ref{k-mean-values}) implies further that  $k^\circ$ is well defined almost everywhere, and the measure $\mu_f$ has the form
\[
d\mu_f=k^\circ\,d\tau+d\sigma,
\]
where $\sigma$ is singular with respect to normalized Lebesgue measure $\tau$. 

(iv) Since $\tilde{f}:\overline{\ddd}^d\rightarrow\ddd$, the function $d$ defined by (\ref{k-real-part}) extends to a bounded, polyharmonic function $\tilde{k}:\overline{\ddd}^d\rightarrow\real_+$. It follows that 
\[
k(z)=\int_{\torus^d}K_z\tilde{k}\,d\tau\qquad(z\in\torus^d)
\]
and hence that $d\mu_f=\tilde{k}\,d\tau$ (see \cite[\S2.3.2(a)]{Ru:1969}); $\tilde{k}$ is continuous on $\torus^d$  by hypothesis. \end{pf}

\begin{defn}\label{defn-measure}
For standard $f\in\mathcal{S}_d$, call the measure $\mu_f$ of Proposition~\ref{prop-poisson-facts}\ref{associated-measure} the probability measure on $\torus^d$ determined by $f$. 
\end{defn}

For standard $f\in\mathcal{S}_d$, the probability measure $\mu_f$ turns out to have an associated sequence of orthogonal polynomials in $d$-variables. 

\section{Orthogonal polynomials on the torus\label{sec-OPOT}}

The goal here is to represent the convergents $f_n$ of a given standard Schur function $f\in\mathcal{S}_d$ in terms of polynomials on $\complex^d$. Recurrence relations among the polynomials facilitate analysis of the $f_n$ and thereby of $f$. 

\subsection{Construction\label{sec-construction}}

Let $\widehat{\complex}=\complex\cup\{\infty\}$ denote the Riemann sphere. Define 
\begin{equation}\label{complex-projection}
\pro:\complex^2\setminus\{0\}\rightarrow\widehat{\complex},\qquad\pro\binom{u}{v}=u/v.
\end{equation}
Given $M\in GL(2,\complex)$, let $\rho_M$ denote the linear fractional transformation
\begin{equation}\label{linear-fractional-transformation}
\rho_M:\complex\rightarrow\widehat{\complex},\qquad\rho_M(\xi)=\pro\left(M\binom{\xi}{1}\right),
\end{equation}
so that matrix multiplication corresponds to composition of transformations, 
\begin{equation}\label{composition-property}
\rho_{M_1M_2}=\rho_{M_1}\circ\rho_{M_2}.  
\end{equation}

Let $d\geq 1$.  Fix a standard Schur function $f\in\mathcal{S}_d$, with Schur parameters $r:\integer_+\rightarrow\ddd$,  variable allocation map $\nu:\integer_{>0}\rightarrow\{1,\ldots,d\}$ and covergents $f_n$ $(n\geq 1)$. Set 
\begin{equation}\label{M0-Mj}
M_0=\begin{pmatrix}1&1\\ -1&1\end{pmatrix}\quad\mbox{ and }\quad M_n=\begin{pmatrix}z_{\nu_n}&r_nz_{\nu_n}\\ \overline{r}_n&1\end{pmatrix}\quad(n\geq 1).
\end{equation}
Observe that for $h\in\mathcal{H}_d$ and $z\in\ddd^d$, 
\begin{equation}\label{matrix-composition-formulation}
\rho_{M_n}\circ h(z)=\mathfrak{m}_{\nu_n}\mathfrak{g}_{r_n}h(z),
\end{equation}
from which it follows by (\ref{convergents}) that
\begin{equation}\label{fn-formula-1}
f_n(z)=\rho_{M_1}\circ\cdots\circ\rho_{M_n}(0)=\pro\left(M_1\cdots M_n\binom{0}{1}\right)\qquad(n\geq 1).
\end{equation}
\begin{defn}\label{defn-oput}
For each $n\geq 0$ denote the entries of the product $P_n:=M_0\cdots M_n$ as 
\begin{equation}\label{oput}
P_n=\begin{pmatrix}\Psi_n(z)&\Psi_n^\ast(z)\\ -\Phi_n(z)&\Phi_n^\ast(z)\end{pmatrix}
\end{equation}
Each of $\Phi_n(z),\Psi_n(z),\Phi_n^\ast(z),\Psi_n^\ast(z)$ is a polynomial in variables $z\in\complex^d$.  Call $\Phi_n$ $(n\geq 1)$ the sequence of orthogonal polynomials determined by  $f$.  
\end{defn}
Remarks. The awkward seeming notation used to label the polynomial entries of $P_n$ is chosen for consistency with the extant literature on orthogonal polynomials \cite{Si:1OPUC2005,Kh:2008,Si:2011}.  Orthogonality of the $\Phi_n$ with respect to $\mu_f$ is proved below in \S\ref{sec-orthogonality}.

Combining (\ref{fn-formula-1}) and (\ref{oput}) yields 
\begin{equation}\label{fn-formula-2}
f_n=\pro\left(M_0^{-1}\begin{pmatrix}\Psi_n&\Psi_n^\ast\\ -\Phi_n&\Phi_n^\ast\end{pmatrix}\binom{0}{1}\right)=\frac{\Psi_n^\ast-\Phi_n^\ast}{\Psi_n^\ast+\Phi_n^\ast}. 
\end{equation}
and
\begin{equation}\label{fn-formula-3}
\frac{1+f_n}{1-f_n}=\rho_{M_0}\circ f_n=\pro\left(\begin{pmatrix}\Psi_n&\Psi_n^\ast\\ -\Phi_n&\Phi_n^\ast\end{pmatrix}\binom{0}{1}\right)=\frac{\Psi_n^\ast}{\Phi_n^\ast}.
\end{equation}
Remark. Although originally defined on $\ddd^d$, the above representations show $f_n$ and $(1+f_n)/(1-f_n)$ extend to rational functions on $\complex^d$; these extensions will henceforth be denoted by the same $f_n$ and $(1+f_n)/(1-f_n)$, without further notational adornments. 

\subsection{Basic properties\label{sec-properties}}

A key property of Definition~\ref{defn-oput} is the recurrence relation \[P_{n+1}=P_nM_{n+1},\] which, in expanded form (writing $\Phi_n$ in place of $\Phi_n(z)$ etc.~for readability) is 
\begin{equation}\label{matrix-recurrence}
\begin{pmatrix}\Psi_{n+1}&\Psi_{n+1}^\ast\\ -\Phi_{n+1}&\Phi_{n+1}^\ast\end{pmatrix}=
\begin{pmatrix}z_{\nu_{n+1}}\Psi_n+\overline{r}_{n+1}\Psi_n^\ast&r_{n+1}z_{\nu_{n+1}}\Psi_n+\Psi_n^\ast\\ 
-z_{\nu_{n+1}}\Phi_n+\overline{r}_{n+1}\Phi_n^\ast&-r_{n+1}z_{\nu_{n+1}}\Phi_n+\Phi_n^\ast\end{pmatrix}\qquad(n\geq 0).
\end{equation}
Combining the $(2,1)$ and $(2,2)$-entries of (\ref{matrix-recurrence}) yields 
\begin{equation}\label{downward-induction}
(1-|r_{n+1}|^2)\Phi_n^\ast=\Phi_{n+1}^\ast+r_{n+1}\Phi_{n+1},
\end{equation}
while the $(1,2)$ and $(2,2)$-entries combine to yield
\begin{equation}\label{upward-induction}
\Psi_{n+1}^\ast+\Phi_{n+1}^\ast=r_{n+1}z_{\nu_{n+1}}\left(\Psi_n-\Phi_n\right)+ \Psi_n^\ast+\Phi_n^\ast.  
\end{equation}
Writing $\det P_n$ in terms of $\det M_n=z_{\nu_n}(1-|r_n|^2)$ $(n\geq1)$ and $\det M_0=2$ yields 
\begin{equation}\label{determinantal-identity}
\Psi_n\Phi_n^\ast+\Phi_n\Psi_n^\ast=2z_{\nu_1}\cdots z_{\nu_n}(1-|r_1|^2)\cdots(1-|r_n|^2)\qquad (n\geq0).
\end{equation}
Next, set 
\begin{equation}\label{Q-D}
Q=\begin{pmatrix}0&1\\ 1&0\end{pmatrix}\quad\mbox{ and }\quad D=\begin{pmatrix}-1&0\\ 0&1\end{pmatrix}
\end{equation}
and note that 
\begin{equation}\label{Q-D-formulas}
M_0D=QM_0\quad\mbox{ and }\quad M_0Q=-DM_0. 
\end{equation}
Replacing $r_n$ with $-r_n$ in $M_n$ yields the matrix 
\begin{equation}\label{Mj-minus}
\widetilde{M}_n:=\begin{pmatrix}z_{\nu_n}&-r_nz_{\nu_n}\\ -\overline{r}_n&1\end{pmatrix}=DM_nD. 
\end{equation}
Therefore 
\[
\widetilde{P}_n:=M_0\widetilde{M}_1\cdots\widetilde{M}_n=QP_nD=\begin{pmatrix}\Phi_n(z)&\Phi_n^\ast(z)\\ -\Psi_n(z)&\Psi_n^\ast(z)\end{pmatrix}.
\]
Thus replacing $r$ with $-r$ has the effect of interchanging each $\Psi_n$ and $\Phi_n$.

Given $z\in\complex^d$, let $1/z$ denote $(1/z_1,\ldots,1/z_d)$.  
With this notation, 
\[
\begin{split}
z_{\nu_1}\cdots z_{\nu_n}\overline{\begin{pmatrix}\Psi_n(1/\overline{z})&\Psi_n^\ast(1/\overline{z})\\ -\Phi_n(1/\overline{z})&\Phi_n^\ast(1/\overline{z})\end{pmatrix}}&=
z_{\nu_1}\cdots z_{\nu_n}M_0\begin{pmatrix}1/z_{\nu_1}&\overline{r}_1/z_{\nu_1}\\ r_1&1\end{pmatrix}\cdots
\begin{pmatrix}1/z_{\nu_n}&\overline{r}_n/z_{\nu_n}\\ r_n&1\end{pmatrix}\\
&=M_0\begin{pmatrix}1&\overline{r}_1\\ r_1z_{\nu_1}&z_{\nu_1}\end{pmatrix}\cdots
\begin{pmatrix}1&\overline{r}_n\\ r_nz_{\nu_n}&z_{\nu_n}\end{pmatrix}\\
&=M_0(QM_1Q)\cdots(QM_nQ)\\
&=-DP_nQ\quad\mbox{ by (\ref{Q-D-formulas})}\\
&=\begin{pmatrix}\Psi_n^\ast(z)&\Psi_n(z)\\ -\Phi_n^\ast(z)&\Phi_n(z)\end{pmatrix}.
\end{split}
\]
Thus
for each $n\geq 1$,
\begin{equation}\label{star-identity}
\Psi_n^\ast(z)=z_{\nu_1}\cdots z_{\nu_n}\overline{\Psi_n(1/\overline{z})}\qquad\mbox{ and }\qquad
\Phi_n^\ast(z)=z_{\nu_1}\cdots z_{\nu_n}\overline{\Phi_n(1/\overline{z})}.
\end{equation}
\begin{prop}\label{prop-zero-free}
For every $n\geq 0$, both $\Phi_n^\ast$ and $\Psi_n^\ast+\Phi_n^\ast$ are zero free on $\overline{\ddd}^d$. 
\end{prop}
\begin{pf}
Consider first $\Phi_n^\ast$. 
Note that for every $n\geq 0$ and $z\in\torus^d$, $|\Phi_n(z)/\Phi_n^\ast(z)|=1$, by (\ref{star-identity}).  Suppose that $\Phi_n^\ast$ is zero-free on $\overline{\ddd}^d$, and, referring to (\ref{matrix-recurrence}), consider 
\[
\Phi_{n+1}^\ast=-r_{n+1}z_{\nu_{n+1}}\Phi_n+\Phi_n^\ast.
\]
If $\Phi_{n+1}^\ast(z)=0$ for some $z\in\overline{\ddd}^d$, then
\begin{equation}\label{equal1}
r_{n+1}z_{\nu_{n+1}}\frac{\Phi_n(z)}{\Phi_n^\ast(z)}=1.
\end{equation}
But $|\Phi_n/\Phi_n^\ast|\leq 1$ on $\overline{\ddd}^d$ by the maximum modulus principle.  Furthermore $|z_{\nu_{n+1}}|\leq 1$ and $|r_{n+1}|<1$, therefore 
\[
\left|r_{n+1}z_{\nu_{n+1}}\frac{\Phi_n(z)}{\Phi_n^\ast(z)}\right|<1,
\]
contradicting (\ref{equal1}).  Thus $\Phi_{n+1}^\ast$ is zero-free on $\overline{\ddd}^d$. Since $\Phi_0^\ast=1$ is zero-free, the desired result follows by induction. 

The justification that $\Psi_n^\ast+\Phi_n^\ast$ is zero-free on $\ddd^d$ is similar.  Recall from (\ref{fn-formula-2}) that
\[
f_n=\frac{\Psi_n^\ast-\Phi_n^\ast}{\Psi_n^\ast+\Phi_n^\ast}. 
\]
It follows easily from the representation (\ref{fn-formula-1}) that $|f_n|<1$ on $\overline{\ddd}^d$, since each $r_j$ has modulus less than 1. Formulas (\ref{star-identity}) imply that for $z\in\torus^d$, 
\[
|\Psi_n^\ast(z)-\Phi_n^\ast(z)|=|\Psi_n(z)-\Phi_n(z)|,
\]
whereby, for $z\in\torus^d$,
\[
\left|\frac{\Psi_n(z)-\Phi_n(z)}{\Psi_n^\ast(z)+\Phi_n^\ast(z)}\right|=|f_n(z)|<1.
\]
It follows by the maximum principle that 
\begin{equation}\label{strictly-less}
\left|\frac{\Psi_n(z)-\Phi_n(z)}{\Psi_n^\ast(z)+\Phi_n^\ast(z)}\right|<1
\end{equation}
for all $z\in\overline{\ddd}^d$.  Therefore by (\ref{upward-induction}), for $z\in\overline{\ddd}^d$, 
\[
\begin{split}
|\Psi_{n+1}^\ast(z)+\Phi_{n+1}^\ast(z)|
&=|r_{n+1}z_{\nu_{n+1}}\left(\Psi_n(z)-\Phi_n(z)\right)+ \Psi_n^\ast(z)+\Phi_n^\ast(z)|\\
&\geq|\Psi_n^\ast(z)+\Phi_n^\ast(z)|-|\Psi_n(z)-\Phi_n(z)|\\
&>0\qquad\mbox{ by (\ref{strictly-less})}. 
\end{split}
\]
Thus $\Psi_{n+1}^\ast+\Phi_{n+1}^\ast$ is zero free on $\overline{\ddd}^d$.  Since $\Psi_0^\ast+\Phi_0^\ast=2$, the desired result follows by induction. 
\end{pf}

The following proposition is based on the simple observation that if 
\[
r_{j+1}=\cdots=r_{j+m-1}=0
\]
then
\[
\begin{split}
M_{j+1}\cdots M_{j+m}&=\begin{pmatrix}z_{j+1}&0\\ 0&1\end{pmatrix}\cdots\begin{pmatrix}z_{j+m-1}&0\\ 0&1\end{pmatrix}\begin{pmatrix}z_{j+m}&r_{j+m}z_{j+m}\\ \overline{r_{j+m}}&1\end{pmatrix}\\
&=\begin{pmatrix}z_{j+1}\cdots z_{j+m}&r_{j+m}z_{j+1}\cdots z_{j+m}\\ \overline{r_{j+m}}&1\end{pmatrix}.
\end{split}
\]

\begin{prop}\label{prop-monomial-substitution}
Fix $\kappa:\{1,\ldots, d\}\rightarrow\integer_+^D\setminus\{0\}$ for some $d,D\geq 1$, and, given $z\in\complex^D$, write 
\[
z^\kappa=(z^{\kappa(1)},\ldots,z^{\kappa(d)})\in\complex^d.
\]
For $h\in\mathcal{H}_d$, define $h^\kappa\in\mathcal{H}_D$ by $h^\kappa(z)=h(z^\kappa)$.  
If $f\in\mathcal{S}_d$, then $f^\kappa\in\mathcal{S}_D$. 
\end{prop} 
\begin{pf}
Let $f\in\mathcal{S}_d$ have Schur parameters $r$ and variable allocation function $\nu$. For $z=(z_1,\ldots,z_N)\in\complex^D$ let $z^{\kappa\circ\nu}$ denote the infinite sequence of multinomials
\[
z^{\kappa\circ\nu}=\bigl(z^{\kappa(\nu(1))},z^{\kappa(\nu(2))},\ldots\bigr).
\]
For each $j\geq 1$, set $m_j=\bigl|\kappa(\nu(j))\bigr|\geq 1$, and fix a mapping 
\[
\sigma_j:\{1,\ldots,m_j\}\rightarrow\{1,\ldots,D\}
\]
such that 
\[
z_{\sigma_j(1)}\cdots z_{\sigma_j(m_j)}=z^{\kappa(\nu(j))}.
\]
(Thus in effect $\sigma_j$ is an elementary ordered partition of the multi-index $\kappa(\nu(j))$.)  Define $\widetilde{r}:\integer_+\rightarrow\overline{\ddd}$ in terms of the given sequence $r$ by replacing each $r_j$ with a block 
\[
(0,\ldots,0,r_j)
\]
consisting of $m_j-1$ zeros followed by $r_j$.  Define 
\[
\widetilde{\nu}:\integer_{>0}\rightarrow\{1,\ldots,D\}
\]
as follows.  To begin, set $n_0=0$, and, for each $j\geq 1$, set $n_j=m_1+\cdots+m_j$.  Then set 
\[
\widetilde{\nu}(n_{j-1}+s)=\sigma_j(s)\qquad(1\leq j; 1\leq s\leq m_j).
\]
Letting $\tilde{f}$ correspond to $\tilde{r}$ and $\tilde{\nu}$, it follows directly from the above construction that for each $j\geq 1$,
\[
\tilde{f}_j(z)=f_j(z^\kappa)\qquad(z\in\ddd^D).
\]
Therefore, setting $\widetilde{f}=\lim_{n\rightarrow\infty}\tilde{f}_n\in\mathcal{S}_D$,
it follows that $f^\kappa=\widetilde{f}\in\mathcal{S}_D$.  
\end{pf}

\subsection{Orthogonality with respect to $\mu_f$\label{sec-orthogonality}}

For each $n\geq 1$, denote by $Z_n^\nu$ the set of monomials 
\begin{equation}\label{monomials}
Z_n^\nu=\left\{z_{\nu_{\eta_1}}\cdots z_{\nu_{\eta_j}}\,\left|\,1\leq j\leq n\;\;\&\;\;1\leq \eta_1<\cdots<\eta_j\leq n\right.\right\}.
\end{equation}
In addition, set $Z_0^\nu=\emptyset$.  Observe that $\spn\left(Z_n^\nu\cup\{1\}\right)$ contains each of the entries of $P_n$, as follows readily from (\ref{matrix-recurrence}). 

Note that because each $f_n$ is holomorphic of modulus strictly less than 1 on $\overline{\ddd}^d$, it determines a probability measure $\mu_{f_n}$ on $\torus^d$ that is absolutely continuous with respect to Lebesgue measure.  And, by Proposition~\ref{prop-real-part-rep}, $d\mu_{f_n}$ has the form 
\begin{equation}\label{mu-fn}
d\mu_{f_n}=\Re\frac{1+f_n(z)}{1-f_n(z)}d\tau
=\frac{\prod_{j=1}^n(1-|r_j|^2)}{|\Phi_n^\ast(z)|^2}d\tau\qquad(z\in\torus^d).
\end{equation}

\begin{prop}\label{prop-monomial-induction}
Let $\rho$ be a measure on $\mathbb{T}^d$. For any $n\geq 1$, if $\Phi_n^\ast\perp Z_n^\nu$ in $L^2(d\rho)$, then $\Phi_n\perp Z_{n-1}^\nu\cup\{1\}$ and $\Phi_{n-1}^\ast\perp Z_{n-1}^\nu$. 
\end{prop}
\begin{pf}
Suppose $\Phi_n^\ast\perp Z^\nu_n$ for some $n\geq 1$. 
Given an arbitrary $p\in Z^\nu_{n-1}\cup\{1\}$, note that $q=z_{\nu_1}\cdots z_{\nu_n}/p\in Z^\nu_n$. Note also that $q=z_{\nu_1}\cdots z_{\nu_n}\overline{p}$ on $\mathbb{T}^d$.  Therefore by (\ref{star-identity})
\[
\int_{\mathbb{T}^d}\Phi_n\overline{p}\,d\rho=\int_{\mathbb{T}^d}\overline{\Phi_n^\ast}z_{\nu_1}\cdots z_{\nu_n}\overline{p}\,d\rho
=\int_{\mathbb{T}^d}\overline{\Phi_n^\ast}q\,d\rho=0,
\]
since $\Phi_n^\ast\perp Z^\nu_n$.  Thus $\Phi_n\perp Z^\nu_{n-1}\cup\{1\}$.  Now, $\Phi_{n-1}^\ast\in\spn\{\Phi_n^\ast,\Phi_n\}$ by (\ref{downward-induction}), and hence
\[
\Phi_{n-1}^\ast\perp Z^\nu_n\cap\left(Z^\nu_{n-1}\cup\{1\}\right)=Z^\nu_{n-1}.
\] 
\end{pf}

\begin{prop}\label{prop-base-orthogonality}
For every $n\geq 1$, $\Phi_n^\ast\perp Z_n^\nu$ in $L^2(d\mu_{f_n})$. 
\end{prop}
\begin{pf}  
Let $p\in Z_n^\nu$, and note that $p(0)=0$ since $p$ is a non-constant monomial.  By (\ref{mu-fn}),
\[
\Phi_n^\ast\overline{p}\,d\mu_{f_n}=\Phi_n^\ast\overline{p}\,\frac{\prod_{j=1}^n(1-|r_j|^2)}{\Phi_n^\ast\overline{\Phi_n^\ast}}\,d\tau=\prod_{j=1}^n(1-|r_j|^2)\frac{\overline{p}}{\overline{\Phi_n^\ast}}\,d\tau.
\]
Since $\overline{p}/\overline{\Phi_n^\ast}$ is antiholomorphic on $\overline{\ddd}^d$ by Proposition~\ref{prop-zero-free}, it follows that 
\[
\int_{\mathbb{T}^d}\Phi_n^\ast\overline{p}\,d\mu_{f_n}=\prod_{j=1}^n(1-|r_j|^2)\int_{\mathbb{T}^d}\frac{\overline{p}}{\overline{\Phi_n^\ast}}\,d\tau=\prod_{j=1}^n(1-|r_j|^2)\frac{\overline{p(0)}}{\overline{\Phi_n^\ast(0)}}=0,
\]
since $\Phi_n^\ast(0)=1$ and $p\in Z_n^\nu$.  Thus $\Phi_n^\ast\perp Z_n^\nu$. 
\end{pf}

Starting from the result of Proposition~\ref{prop-base-orthogonality}, repeated application of Proposition~\ref{prop-monomial-induction} yields the following. 
\begin{cor}\label{cor-orthogonality}
For every $1\leq j\leq n$, $\Phi_j^\ast\perp Z^\nu_j$ and $\Phi_j\perp Z^\nu_{j-1}\cup\{1\}$ in $L^2(d\mu_{f_n})$. 
\end{cor}
\begin{prop}\label{prop-phi-star}
For every $n\geq 1$ and $1\leq j\leq n$, \[\int_{\mathbb{T}^d}\Phi_j^\ast\,d\mu_{f_n}=\prod_{s=1}^j(1-|r_s|^2).\]
\end{prop}
\begin{pf}
Fix $n\geq 1$. Given $1\leq j\leq n$, equation (\ref{downward-induction}) yields
\[
\int_{\mathbb{T}^d}\Phi_{j-1}^\ast\,d\mu_{f_n}=\frac{1}{1-|r_j|^2}\int_{\mathbb{T}^d}\Phi_j^\ast+r_j\Phi_j\,d\mu_{f_n}=
\frac{1}{1-|r_j|^2}\int_{\mathbb{T}^d}\Phi_j^\ast\,d\mu_{f_n},
\]
the latter equality since $\Phi_j\perp 1$, by Corollary~\ref{cor-orthogonality}. Therefore
\[
\int_{\mathbb{T}^d}\Phi_{j-1}^\ast\,d\mu_{f_n}=\frac{1}{\prod_{s=j}^n(1-|r_s|^2)}\int_{\mathbb{T}^d}\Phi_n^\ast\,d\mu_{f_n}.
\]
But by the formula~(\ref{mu-fn}),
\[
\int_{\mathbb{T}^d}\Phi_n^\ast\,d\mu_{f_n}=\int_{\mathbb{T}^d}\Phi_n^\ast\,\frac{\prod_{j=1}^n(1-|r_j|^2)}{\Phi_n^\ast\overline{\Phi_n^\ast}}\,d\tau=\frac{\prod_{j=1}^n(1-|r_j|^2)}{\overline{\Phi_n^\ast(0)}}=\prod_{j=1}^n(1-|r_j|^2)
\]
since $1/\overline{\Phi_n^\ast}$ is antiholomorphic on $\overline{\ddd}^d$ and $\Phi_n^\ast(0)=1$.  The desired result follows.
\end{pf}

\begin{prop}\label{prop-orthogonality}
Let $n\geq 1$ and $j,k\in\{0,\ldots,n\}$. Then 
\[
\int_{\mathbb{T}^d}\Phi_j\overline{\Phi_k}\,d\mu_{f_n}=\left\{
\begin{array}{cc}
0&\mbox{ if }j\neq k\\
\prod_{s=1}^j(1-|r_s|^2)&\mbox{ if }j=k
\end{array}\right..
\]
\end{prop}
\begin{pf}
Note that for each $j$ in the range $1\leq j\leq n$, $\overline{\Phi_j^\ast}=1+\overline{p_j}$ for some $p_j\in \spn Z^\nu_j$, since $\Phi_j^\ast\in\spn\left(Z^\nu_j\cup\{1\}\right)$ and $\Phi_j^\ast(0)=1$. By Corollary~\ref{cor-orthogonality} it follows that 
\[
\int_{\mathbb{T}^d}\Phi_j^\ast\overline{\Phi_j^\ast}\,d\mu_{f_n}=\int_{\mathbb{T}^d}\Phi_j^\ast\,d\mu_{f_n}=\prod_{s=1}^j(1-|r_s|^2),
\]
the latter equality by Proposition~\ref{prop-phi-star}.  Therefore 
\[
\int_{\mathbb{T}^d}\Phi_j\overline{\Phi_j}\,d\mu_{f_n}=\int_{\mathbb{T}^d}\Phi_j^\ast\overline{\Phi_j^\ast}\,d\mu_{f_n}=\prod_{s=1}^j(1-|r_s|^2)
\]
since $\Phi_j\overline{\Phi_j}=\Phi_j^\ast\overline{\Phi_j^\ast}$ on $\mathbb{T}^d$.  On the other hand, given any distinct $j,k\in\{0,\ldots,n\}$, suppose without loss of generality that $j>k$.  Corollary~\ref{cor-orthogonality} implies $\Phi_j\perp\Phi_k$ in $L^2(d\mu_{f_n})$ since $\Phi_k\in\spn\left(Z^\nu_{j-1}\cup\{1\}\right)$, which proves
\[
\int_{\mathbb{T}^d}\Phi_j\overline{\Phi_k}\,d\mu_{f_n}=0.
\]
\end{pf}

Thus the sequence of polynomials $\Phi_0,\ldots,\Phi_n$ is orthogonal on $\mathbb{T}^d$ with respect to $d\mu_{f_n}$.  Note however that $\spn\{\Phi_0,\ldots,\Phi_n\}$ may be a proper subspace of $\spn\left(Z^\nu_n\cup\{1\}\right)$. For example, if $d\geq n\geq 2$ and the restriction of $\nu$ to $\{1,\ldots,n\}$ is injective, then
\[
\dim\spn\left(Z^\nu_n\cup\{1\}\right)=2^n>n+1=\dim\spn\{\Phi_0,\ldots,\Phi_n\}.
\]
So although $\Phi_0,\ldots,\Phi_n$ are orthogonal on $\mathbb{T}^d$ with respect to $d\mu_{f_n}$, they do not in general span the space of polynomials in $d$ variables of degree at most $n$. 

\begin{thm}\label{thm-orthogonality}
The sequence of measures $\mu_{f_n}$ $(n\geq 1)$ converges weakly to $\mu_f$, and the polynomials 
$
\Phi_n$ $(n\geq 0)
$
are orthogonal on $\mathbb{T}^d$ with respect to $\mu_f$.  More precisely, for any integers $j,k\geq 0$, 
\[
\int_{\mathbb{T}^d}\Phi_j\overline{\Phi_k}\,d\mu_f=\left\{
\begin{array}{cc}
0&\mbox{ if }j\neq k\\
\prod_{s=1}^j(1-|r_s|^2)&\mbox{ if }j=k
\end{array}\right..
\]
\end{thm}
\begin{pf}
By Alaoglu's theorem, the sequence of probability measures $\mu_{f_n}$ has a subsequence $\mu_{f_{j_n}}$ $(n\geq1)$ converging weakly in the dual space to $C(\torus^d)$ to some limit $\mu_\ast$. For each fixed 
\[
z=(s_1e^{i\varphi_1},\ldots,s_de^{i\varphi_d})\in\ddd^d,
\]
the corresponding Poisson kernel $K_z$ given by (\ref{poisson-kernel})
belongs to $C(\torus^d)$.  Also, Theorem~\ref{thm-full-taylor} implies
\[
\Re\frac{1+f(z)}{1-f(z)}=\lim_{n\rightarrow\infty}\Re\frac{1+f_{j_n}(z)}{1-f_{j_n}(z)}. 
\]
Therefore 
\[
\int_{\torus^d}K_z\,d\mu_\ast=\lim_{n\rightarrow\infty}\int_{\torus^d}K_z\,d\mu_{f_{j_n}}
=\lim_{n\rightarrow\infty}\Re\frac{1+f_{j_n}(z)}{1-f_{j_n}(z)}=\Re\frac{1+f(z)}{1-f(z)}=\int_{\torus^d}K_z\,d\mu_f.
\]
By uniqueness of boundary measures \cite[Thm.~2.1.3(e)]{Ru:1969}, it follows that $\mu_\ast=\mu_f$, and moreover that $\mu_f$ is the unique limit point of the sequence $\mu_{f_n}$. Thus $\mu_{f_n}\rightarrow\mu_f$ weakly.  In particular,
for any non-negative integers $j,k$, one has $\Phi_j\overline{\Phi_k}\in C(\torus^d)$, and so
\[
\int_{\mathbb{T}^d}\Phi_j\overline{\Phi_k}\,d\mu_f=\lim_{n\rightarrow\infty}\int_{\mathbb{T}^d}\Phi_j\overline{\Phi_k}\,d\mu_{f_n}=\left\{
\begin{array}{cc}
0&\mbox{ if }j\neq k\\
\prod_{s=1}^j(1-|r_s|^2)&\mbox{ if }j=k
\end{array}\right.,
\]
by Proposition~\ref{prop-orthogonality}. 
\end{pf}

Remark. In the case $d=1$ this proves a celebrated theorem of Geronimus \cite{Ge:1944} identifying Schur parameters with recurrence coefficients. See \cite[Thm.~8.20]{Kh:2008}, \cite[\S3.2]{Si:1OPUC2005} and \cite[Thm.~2.5.2]{Si:2011}. 

\subsection{Further results\label{sec-further}}

\begin{prop}\label{prop-real-part-rep}
For each $n\geq 0$, for every $z\in\torus^d$, 
\begin{equation}\label{real-part-rep}
\Re\frac{1+f_n(z)}{1-f_n(z)}=\frac{\prod_{j=1}^n(1-|r_j|^2)}{|\Phi_n^\ast(z)|^2}=\frac{\prod_{j=1}^n(1-|r_j|^2)}{|\Phi_n(z)|^2}.
\end{equation}
\end{prop}
\begin{pf} This follows directly from (\ref{fn-formula-3}), (\ref{determinantal-identity}) and (\ref{star-identity}).\end{pf}
\begin{prop}\label{prop-transmission}
For each $n\geq 1$, set 
\begin{equation}\label{gn}
g_n=\frac{2\prod_{j=1}^n\sqrt{1-|r_j|^2}}{\Psi_n^\ast+\Phi_n^\ast}.
\end{equation}
For every $z\in\torus^d$, 
\begin{equation}\label{transmission}
1-|f_n(z)|^2=|g_n(z)|^2.
\end{equation}
\end{prop}
\begin{pf} Formulas (\ref{fn-formula-2}), (\ref{star-identity}) and (\ref{determinantal-identity}) yield the stated result.\end{pf}

By Proposition~\ref{prop-zero-free}, $g_n$ is holomorphic and zero free on $\overline{\ddd}^d$. Applying the mean value theorem to the harmonic function $\log|g_n|$ therefore yields the following. 
\begin{cor}\label{cor-fn-integral}
\[
\int_{\torus^d}\log(1-|f_n|^2)\,d\tau=\log\prod_{j=1}^n(1-|r_j|^2).
\]
\end{cor}

\begin{prop}\label{prop-fn-kn-formula}
Given $n\geq 0$, let $f_n$ and $k_n$ be the corresponding convergent and fractional part of $f$, as per Definition~\ref{defn-continued-fraction}, and set
\[
B=\frac{\Psi_n-\Phi_n}{\Psi_n^\ast+\Phi_n^\ast}.
\]
Then, on the torus $\torus^d$,  $|B|=|f_n|$, and 
\begin{equation}\label{fn-kn-formula}
\bigl(1-|f^\circ|^2\bigr)\;\bigl|1+Bk_n^\circ\bigr|^2=\bigl(1-|f_n|^2\bigr)\bigl(1-|k_n^\circ|^2\bigr). 
\end{equation}
\end{prop}
\begin{pf} Consider the factorization $f=\mathfrak{f}_nk_n$ from (\ref{f-factorizations}), noting that $f_n=\mathfrak{f}_n\mathbf{c}_0$. It follows from (\ref{matrix-composition-formulation}) and the definition (\ref{frak-fn}) of $\mathfrak{f}_n$ that 
\[
\begin{split}
\mathfrak{f}_nk_n&=\pro\left(M_1\cdots M_n\binom{k_n}{1}\right)\\
&=\pro\left(M_0^{-1}\begin{pmatrix}\Psi_n&\Psi_n^\ast\\ -\Phi_n&\Phi_n^\ast\end{pmatrix}\binom{k_n}{1}\right)\\
&=\frac{(\Psi_n+\Phi_n)k_n+\Psi_n^\ast-\Phi_n^\ast}{(\Psi_n-\Phi_n)k_n+\Psi_n^\ast+\Phi_n^\ast}\\
&=C\frac{A+k_n}{1+Bk_n},
\end{split}
\]
where 
\[
A=\frac{\Psi_n^\ast-\Phi_n^\ast}{\Psi_n+\Phi_n},\quad B=\frac{\Psi_n-\Phi_n}{\Psi_n^\ast+\Phi_n^\ast}\quad\mbox{ and }\quad C=\frac{\Psi_n+\Phi_n}{\Psi_n^\ast+\Phi_n^\ast}. 
\]
Formula (\ref{fn-formula-2}) and the identities (\ref{star-identity}) imply that for $z\in\torus^d$, $|C(z)|=1$, 
\[
|A(z)|=|B(z)|=|f_n(z)|, 
\]
and $B(z)=\overline{A(z)}$.  It follows that on $\torus^d$,
\[
1-|f^\circ|^2=1-\left|\frac{A+k_n^\circ}{1+Bk_n^\circ}\right|^2=\frac{\bigl(1-|f_n|^2\bigr)\bigl(1-|k_n^\circ|^2\bigr)}{\bigl|1+Bk_n^\circ\bigr|^2},
\]
completing the proof. 
\end{pf}

\section{Szeg\H{o}-Verblunsky theorems\label{sec-verblunsky}}

\subsection{Preliminaries\label{sec-preliminaries}}

\begin{lem}\label{lem-general-outer}
Let $h\in\mathcal{H}_d$ satisfy $h(0)=0$. Then
\[
\int_{\torus^d}\log|1+h^\circ|\,d\tau=0.
\]
\end{lem}
\begin{pf}
Any bounded holomorphic function $k:\ddd^d\rightarrow\complex$ having positive real part is outer (see \cite[Thm.~4.4.9, p.~77]{Ru:1969}), meaning 
\begin{equation}\label{general-outer}
\int_{\torus^d}\log|k^\circ|\,d\tau=\log|k(0)|.
\end{equation}
The condition $h(0)=0$ implies $|h(z)|<1$ for all $z\in\ddd^d$; hence the real part of $1+h$ is positive.  Setting $k=1+h$ and applying (\ref{general-outer}) yields the conclusion of the lemma. 
\end{pf}

\begin{cor}\label{cor-easy-outer}
Let $f\in\mathcal{S}_d$ be a standard Schur function.  Then 
\[
\int_{\torus^d}\log|1-f^\circ|\,d\tau=0.
\]
\end{cor}

\begin{prop}\label{prop-f-fn-comparison}
For every $n\geq 1$, 
\[
\int_{\torus^d}\log(1-|f^\circ|^2)d\tau=\log\prod_{j=1}^n(1-|r_j|^2)+\int_{\torus^d}\log(1-|k_n^\circ|^2)\,d\tau.
\]
\end{prop}
\begin{pf}
By Proposition~\ref{prop-fn-kn-formula},
\begin{multline}\nonumber
\int_{\torus^d}\log(1-|f^\circ|^2)d\tau+\int_{\torus^d}\log\left|1+Bk_n^\circ\right|\,d\tau\\
=\int_{\torus^d}\log(1-|f_n|^2)\,d\tau+\int_{\torus^d}\log(1-|k_n^\circ|^2)\,d\tau.
\end{multline}
Since $|B|=|f_n|<1$ on $\torus^d$ and $k_n^\circ(0)=0$, Lemma~\ref{lem-general-outer} implies
\[
\int_{\torus^d}\log\left|1+Bk_n^\circ\right|\,d\tau=0.
\]
The desired result then follows from Corollary~\ref{cor-fn-integral}. 
\end{pf}

\begin{prop}\label{prop-radial-increase}
Let $h\in\mathcal{H}_d$ be non-constant, and set
\begin{equation}\label{Lambda-notation}
\Lambda_h(\varepsilon)=-\int_{\torus^d}\log(1-|h(\varepsilon z)|^2)\,d\tau\qquad(0\leq \varepsilon\leq 1),
\end{equation}
with the understanding that $h$ is to be replaced by $h^\circ$ when $\varepsilon=1$.  The function 
\[
\Lambda_h:[0,1]\rightarrow\real_+\cup\{\infty\}
\]
is strictly increasing on $[0,1]$, and, if $0<\varepsilon<1$, then $0<\Lambda_h(\varepsilon)<\infty$.
\end{prop}
\begin{pf}
For any integer $m\geq 1$, boundedness of $h^m$ implies its radial limit belongs to $L^2(\torus^d)$. Therefore $h^{m\circ}$ has a Fourier series on the torus of the form
\begin{equation}\label{fourier-series}
h^{m\circ}(z)=\sum_{\alpha\in\integer_+^d}a_\alpha z^\alpha\qquad(z\in\torus^d)
\end{equation}
for which Plancherel's theorem asserts
\begin{equation}\label{plancherel}
\int_{\torus^d}\bigl|h^{m\circ}(z)\bigr|^2\,d\tau=\sum_{\alpha\in\integer_+^d}|a_\alpha|^2<\infty. 
\end{equation}
The Fourier expansion (\ref{fourier-series}) is the radial limit of the power series for $h^m$.  Thus, for $0<\varepsilon<1$, 
\begin{equation}\label{fourier-series-2}
h^{m}(\varepsilon z)=\sum_{\alpha\in\integer_+^d}a_\alpha\varepsilon^{|\alpha|}z^\alpha\qquad(z\in\torus^d),
\end{equation}
and Plancherel's theorem yields
\begin{equation}\label{plancherel-2}
\int_{\torus^d}\bigl|h^{m}(\varepsilon z)\bigr|^2\,d\tau=\sum_{\alpha\in\integer_+^d}|a_\alpha\varepsilon^{|\alpha|}|^2<\infty.
\end{equation}
Since $h$ is not identically 0 (being non-constant) at least one Fourier coefficient $a_\alpha$ is non-zero. Therefore the latter series---and hence the integral to which it is equal---is strictly increasing as a function of $\varepsilon\in[0,1]$. It follows from the expansion
\begin{equation}\label{log-expansion}
-\log(1-x)=x+\frac{x^2}{2}+\frac{x^3}{3}+\cdots
\end{equation}
that 
\[
\Lambda_h(\varepsilon)=\sum_{m=1}^\infty \frac{1}{m(2\pi)^d}\int_{\torus^d}\bigl|h^{m}(\varepsilon z)\bigr|^2\,d\tau.
\]
Each term in the expansion is strictly increasing in $\varepsilon$, therefore so is $\Lambda_n(\varepsilon)$.  

On the other hand, if $0<\varepsilon<1$, then $|h(\varepsilon z)|$ is bounded away from 1 for $z\in\torus^d$. Thus 
\[
-\log(1-|h(\varepsilon z)|^2)
\]
is non-negative and bounded, and its integral is finite (and strictly positive by (\ref{plancherel-2})). 
\end{pf}

\subsection{A multivariate Szeg\H{o}-Verblunsky theorem\label{sec-multivariate}}

\begin{thm}\label{thm-boyd}
Let $f$ be a standard Schur function, and $r:\integer_+\rightarrow\ddd$ its sequence of Schur parameters. Then 
\begin{equation}\label{boyd}
\int_{\torus^d}\log(1-|f^\circ|^2)\,d\tau=\sum_{j=1}^\infty\log(1-|r_j|^2).
\end{equation}
Equation (\ref{boyd}) is valid irrespective of whether the left and right-hand sides are finite or infinite.  The values are finite if and only if $r\in\ell^2(\integer_+)$.
\end{thm}
\begin{pf}
Since $-\log(1-x)\geq 0$ if $0\leq x\leq 1$, 
Proposition~\ref{prop-f-fn-comparison} implies 
\begin{equation}\label{lambda-lower-bound}
-\int_{\torus^d}\log(1-|f^\circ|^2)\,d\tau\geq-\sum_{j=1}^n\log(1-|r_j|^2)\qquad(n\geq 1). 
\end{equation}
By a standard argument, $-\sum_{j=1}^\infty\log(1-|r_j|^2)<\infty$ if and only if $r\in\ell^2(\integer_+)$.  In particular, if 
$-\sum_{j=1}^\infty\log(1-|r_j|^2)$ diverges (equivalently, if $\prod_{j=1}^\infty(1-|r_j|^2)=0$), then (\ref{lambda-lower-bound}) forces 
\[
-\int_{\torus^d}\log(1-|f^\circ|^2)\,d\tau=\infty.
\]

Suppose now that $r\in\ell^2(\integer_+)$, whence
\[
\gamma:=-\sum_{j=1}^\infty\log(1-|r_j|^2)<\infty.
\]
In this case (\ref{lambda-lower-bound}) guarantees $-\int_{\torus^d}\log(1-|f^\circ|^2)\,d\tau\geq\gamma$. 
It remains to obtain the reverse inequality. 
Fatou's lemma concerning positive, measurable functions asserts
\begin{equation}\label{fatou-lemma}
-\int_{\torus^d}\log(1-|f^\circ|^2)\,d\tau\leq\liminf_{\varepsilon\rightarrow 1^-}\Lambda_f(\varepsilon)=\lim_{\varepsilon\rightarrow 1^-}\Lambda_f(\varepsilon),
\end{equation}
the latter equality by Proposition~\ref{prop-radial-increase}.  Recall $f_n\rightarrow f$ uniformly on compact sets in $\ddd^d$.  Therefore, for any fixed $0\leq\varepsilon<1$, 
\begin{equation}\label{lambda-limit}
\Lambda_f(\varepsilon)=\lim_{n\rightarrow\infty}\Lambda_{f_n}(\varepsilon).
\end{equation}
On the other hand, Proposition~\ref{prop-radial-increase} and Corollary~\ref{cor-fn-integral} imply that for each fixed $n\geq 1$,
\[
\Lambda_{f_n}(\varepsilon)\leq\Lambda_{f_n}(1)=-\log\prod_{j=1}^n(1-|r_j|^2)\leq\gamma.
\]
Therefore (\ref{lambda-limit}) yields $\Lambda_f(\varepsilon)\leq\gamma$, whereby (\ref{fatou-lemma}) yields $\Lambda_f(1)\leq\gamma$. Thus $\Lambda_f(1)=\gamma$. 
\end{pf}

The multivariate analogue of the classical Szeg\H{o}-Verblunsky theorem is an easy consequence of the foregoing result. 
\begin{thm}[Multivariate Szeg\H{o}-Verblunsky theorem]\label{thm-multivariate-szego}
Let $f\in\mathcal{S}_d$ be a standard Schur function, $r:\integer_+\rightarrow\ddd$ its sequence of Schur parameters, and $\mu_f$ its probability measure. Write 
\[
d\mu_f=w\,d\tau+d\sigma,
\]
where $\sigma$ is singular with respect to normalized Lebesgue measure $\tau$.  Then 
\begin{equation}\label{multivariate-szego-conclusion}
\int_{\torus^d}\log w\,d\tau=\sum_{j=1}^\infty\log(1-|r_j|^2).
\end{equation}
\end{thm}
\begin{pf}
By Proposition~\ref{prop-poisson-facts}\ref{associated-measure} the absolutely continuous part of $\mu_f$ has the form
\begin{equation}\label{w-form}
w=\Re\frac{1+f^\circ}{1-f^\circ}=\frac{1-|f^\circ|^2}{|1-f^\circ|^2}. 
\end{equation}
Therefore 
\[
\int_{\torus^d}\log w\,d\tau=\int_{\torus^d}\log(1-|f^\circ|^2)\,d\tau-2\int_{\torus^d}\log|1-f^\circ|\,d\tau.
\]
Theorem~\ref{thm-boyd} and Corollary~\ref{cor-easy-outer} then give the desired result. 
\end{pf}

\subsection{An almost periodic Szeg\H{o}-Verblunsky theorem\label{sec-almost-periodic}}

To begin, fix 
\[
\eta=(\eta_1,\ldots,\eta_d)\in\real_{>0}^d
\]
with strictly positive entries. 
Set 
\[
D=\dim\spn_\rat\{\eta_1,\ldots,\eta_d\},
\]
the dimension of $\spn_\rat\{\eta_1,\ldots,\eta_d\}$ as a vector space over the rational numbers $\rat$.  Choose a point $q=(q_1,\ldots,q_D)\in\real_+^D$ such that 
\[
\spn_\rat\{q_1,\ldots,q_D\}=\spn_\rat\{\eta_1,\ldots,\eta_d\}.
\]
and such that each $\eta_j$ is a non-negative integer combination of $q_1,\ldots,q_D$.  The following proposition guarantees that this is always possible.

\begin{prop}\label{prop-positive-integer}
Let $\eta\in\real_+^d$ have strictly positive entries, and let 
$
D=\dim\spn_\rat\{\eta_1,\ldots,\eta_d\}.
$
Then there exist $q\in\real_+^D$ and $\alpha_\nu=(\alpha_\nu^1,\ldots,\alpha_\nu^D)\in\integer_+^D$ $(1\leq \nu\leq d)$ such that 
\[
\eta_\nu=\langle\alpha_\nu,q\rangle\qquad(1\leq \nu\leq d).
\]
\end{prop}
\begin{pf}
Choose a positive basis $\{b_1,\ldots,b_D\}$ for $\spn_\rat\{\eta_1,\ldots,\eta_d\}$ (if $\{b_1,\ldots,b_D\}$ is any basis, then so is $\{|b_1|,\ldots,|b_D|\}$), and set $b=(b_1,\ldots,b_D)^t$, viewing $b$ as a column vector.  Let $B\in\rat_{d\times D}$ denote the unique rational $d\times D$ matrix such that 
\[
\eta=Bb. 
\]
Let $q^{(j)}\in\rat^D_+$ $(j\in\integer_+)$ be a sequence of rational (column) vectors converging to $b$, and set 
\[
P_j=\frac{1}{\langle q^{(j)},q^{(j)}\rangle} q^{(j)}q^{(j)\,t}\qquad (j\in\integer_+),
\]
so that $P_j$ is a rational, symmetric $D\times D$ projection matrix.   For $t\in\rat_+$, set 
\[
Q_{j,t}=I+tP_j\qquad(j\in\integer_+).
\]
Setting
\[
\varepsilon=\frac{1}{2\langle b,b\rangle}\min_{\nu,\lambda}\eta_\nu b_\lambda>0,
\]
observe that, by construction, 
\[
BP_j\rightarrow\frac{1}{\langle b,b\rangle}\eta b^t\quad\mbox{ as }\quad j\rightarrow\infty,
\]
and so there is an index $m^\prime$ such that for all $j\geq m^\prime$, each entry of $BP_j$ exceeds $\varepsilon$.  Fix $t\in\rat_+$ sufficiently large that 
\[
BQ_{j,t}=B+tBP_j
\]
has strictly positive entries for all $j\geq m^\prime$. Note that 
\[
Q_{j,t}^{-1}=I-\frac{t}{1+t}P_j,
\]
and hence that 
\[
Q_{j,t}^{-1}b\rightarrow\frac{1}{1+t}b>0\quad\mbox{ as }\quad j\rightarrow\infty.
\]
Fix $m>m^\prime$ sufficiently large that each of the entries of $Q_{m,t}^{-1}b$ is stricly positive.  Finally, since $BQ_{m,t}$ has rational entries, there is a positive integer $s$ such that $sBQ_{m,t}$ has integer entries.  Set 
\[
A=sBQ_{m,t}\quad\mbox{ and }\quad q=s^{-1}Q_{m,t}^{-1}b.
\]
Then $A$ has positive integer entries, $q\in\real_+^D$, and 
\[
\eta=Aq.
\]
Setting $\alpha_\nu$ to be the $\nu$th row of $A$ $(1\leq \nu\leq d)$ completes the proof.\end{pf}

Let $\alpha=(\alpha_1,\ldots,\alpha_d)$ be the sequence of multi-indices whose existence is assured by Proposition~\ref{prop-positive-integer}, and define
\begin{equation}\label{tau-alpha}
\sigma_\alpha:\complex^D\rightarrow\complex^d,\qquad \sigma_\alpha z=(z^{\alpha_1},\ldots,z^{\alpha_d}).
\end{equation}
Note that 
\begin{equation}\label{A-preserves-torus}
\sigma_\alpha(\torus^D)\subset \torus^d\quad\mbox{ and }\quad \sigma_\alpha(\overline{\ddd}^D)\subset\overline{\ddd}^d.
\end{equation}  
Define torus lines
\[
\begin{split}
\ell_q&:\real\rightarrow\torus^D,\quad \ell_q(\omega)=(e^{iq_1\omega},\ldots,e^{iq_D\omega})\quad\mbox{ and }\\
\ell_\eta&:\real\rightarrow\torus^d,\quad\ell_\eta(\omega)=(e^{i\eta_1\omega},\ldots,e^{i\eta_d\omega}),
\end{split}
\]
and observe that 
\begin{equation}\label{Delta-pull-back}
\ell_\eta=\sigma_\alpha\circ\ell_q.
\end{equation}
Note that the line $\ell_q$ is dense on $\torus^D$ since $q_1,\ldots,q_D$ are linearly independent over $\rat$, whereas if $D<d$, the line $\ell_\eta$ is not dense on $\torus^d$. 

Suppose $f\in\mathcal{S}_d$ has a holomorphic extension $\tilde{f}:\overline{\ddd}^d\rightarrow\ddd$, and let $r:\integer_+\rightarrow\ddd$ and $\nu:\integer_{>0}\rightarrow\{1,\ldots,d\}$ be the respective Schur parameters and variable allocation map of $f$. Given $\eta\in\real_{>0}^d$, interpret $\eta\circ\nu:\integer_{>0}\rightarrow\real_{>0}$ in the natural way by the formula 
\[
\eta\circ\nu(j)=\eta_{\nu_j}\qquad(j\in\integer_{>0}). 
\] 
Holomorphy on $\overline{\ddd}^d$ implies the Taylor series for $f$ given in Theorem~\ref{thm-full-taylor} converges absolutely and uniformly on $\overline{\ddd}^d$. Therefore the function
\[
\tilde{f}\circ\ell_\eta:\real\rightarrow\ddd
\]
is represented by the absolutely convergent series
\begin{equation}\label{almost-periodic-series}
\tilde{f}\circ\ell_\eta(\omega)=\sum_{\alpha\in A}c_\alpha(r)e^{i\langle \alpha,\eta\circ\nu\rangle\omega},
\end{equation} 
where 
\[
\langle \alpha,\eta\circ\nu\rangle=\sum_{j=1}^\infty\alpha_j\eta_{\nu_j}. 
\]
Now, the function (\ref{almost-periodic-series}) is almost periodic in the sense of Besicovitch \cite[Ch.~II]{Be:1955}, with almost periods $\langle \alpha,\eta\circ\nu\rangle$ $(\alpha\in A)$.  The purpose of the present section is to prove the analogue of Theorem~\ref{thm-boyd} for such functions.  The following ergodic lemma is a crucial ingredient. 

\begin{lem}\label{lem-continuous-birkhoff}
Let $g\in C(\torus^d)$, and let $\eta\in\real_{>0}^d$ satisfy $\dim_{\rat}\{\eta_1,\ldots,\eta_d\}=d$. Then
\[
\lim_{L\rightarrow\infty}\frac{1}{2L}\int_{-L}^Lg\circ\ell_\eta(\omega)\,d\omega=\int_{\torus^d}g\,d\tau.
\]
\end{lem}
\begin{pf}
Denote $\mathbb{1}=(1,\ldots,1)\in\torus^d$, and let $\varepsilon>0$ be arbitrary. Fix $\delta>0$ such that 
\begin{equation}\label{3-epsilon-0}
\|z-w\|_\infty<\delta\|\eta\|_{\infty}\Rightarrow|g(z)-g(w)|<\varepsilon/3\qquad (z,w\in\torus^d). 
\end{equation}
(Such a $\delta$ exists by continuity of $g$ and compactness of $\torus^d$.)
For any $m\geq 1$ and $z=(e^{i\theta_1},\ldots,e^{i\theta_d})\in\torus^d$, set 
\begin{equation}\label{A-m-delta}
A^\delta_mg(z)=\frac{1}{2m}\sum_{j=-m}^{m-1}g(e^{i(\theta_1+j\delta\eta_1)},\ldots,e^{i(\theta_d+j\delta\eta_d)}).
\end{equation}
Translation on the torus by $\pm\delta\eta$ is measure preserving and has no nontrivial invariant subsets.  Therefore by Birkhoff's ergodic theorem \cite[\S1.2]{Kr:1985},
\begin{equation}\label{birkhoff-limit}
\lim_{m\rightarrow\infty}A^\delta_mg(z)=\int_{\torus^d}g\,d\tau 
\end{equation}
for almost every $z\in\torus^d$. Choose $z\in\torus^d$ such that $\|z-\mathbb{1}\|_\infty<\delta\|\eta\|_\infty$ and (\ref{birkhoff-limit}) holds. Observe that
\begin{equation}\label{3-epsilon-1}
\left|A^\delta_mg(\mathbb{1})-A^\delta_mg(z)\right|<\varepsilon/3
\end{equation}
independently of $m\geq 1$.  Using (\ref{birkhoff-limit}), fix $n$ large enough that
\begin{equation}\label{3-epsilon-2}
\left|A^\delta_mg(z)-\int_{\torus^d}g\,d\tau\right|<\varepsilon/3\qquad(m\geq n). 
\end{equation}
Note that $A^\delta_mg(\mathbb{1})$ is a Riemann sum approximation to the integral 
\[
\frac{1}{2L}\int_{-L}^Lg\circ\ell_\eta(\omega)\,d\omega\quad\mbox{ where }\quad L=m\delta. 
\]
By (\ref{3-epsilon-0}), the error in this approximation satisfies
\begin{equation}\label{3-epsilon-3}
\left|A^\delta_mg(\mathbb{1})-\frac{1}{2L}\int_{-L}^Lg\circ\ell_\eta(\omega)\,d\omega\right|<\varepsilon/3\qquad(m\geq1; L=m\delta). 
\end{equation}
Combining (\ref{3-epsilon-3},\ref{3-epsilon-1},\ref{3-epsilon-2}), for every $m\geq n$ and $L=m\delta$, 
\begin{multline}\label{3-epsilon-final}
\left|\frac{1}{2L}\int_{-L}^Lg\circ\ell_\eta(\omega)\,d\omega-\int_{\torus^d}g\,d\tau\right|\leq\\
\left|A^\delta_mg(\mathbb{1})-\frac{1}{2L}\int_{-L}^Lg\circ\ell_\eta(\omega)\,d\omega\right|+
\left|A^\delta_mg(\mathbb{1})-A^\delta_mg(z)\right|+
\left|A^\delta_mg(z)-\int_{\torus^d}g\,d\tau\right|<\varepsilon. 
\end{multline}
Since (\ref{3-epsilon-final}) holds for all $L=m\delta$ with $m\geq n$, it follows by boundedness of $g$, that
\[
\limsup_{L\rightarrow\infty}\left|\frac{1}{2L}\int_{-L}^Lg\circ\ell_\eta(\omega)\,d\omega-\int_{\torus^d}g\,d\tau\right|\leq\varepsilon. 
\]
Since $\varepsilon$ was arbitrary, the conclusion of the lemma follows. 
\end{pf}

Here is the main almost periodic result.  

\begin{thm}[Almost periodic Szeg\H{o}-Verblunsky theorem]\label{thm-almost-periodic-main}
Let $f\in\mathcal{S}_d$ be a standard Schur function that extends to a holomorphic function $\tilde{f}:\overline{\ddd}^d\rightarrow\ddd$, and let $r:\integer_+\rightarrow\ddd$ be its sequence of Schur parameters.
Then, for every $\eta\in\real_{>0}^d$, 
the function
$
\tilde{f}\circ\ell_\eta:\real\rightarrow\ddd
$
is almost periodic, represented by the absolutely convergent series
(\ref{almost-periodic-series}), and
\begin{equation}\label{main-almost-periodic-result}
\lim_{L\rightarrow\infty}\frac{1}{2L}\int_{-L}^L-\log\left(1-\bigl|\tilde{f}\circ\ell_\eta(\omega)\bigr|^2\right)d\omega=\sum_{j=1}^\infty-\log(1-|r_j|^2).  
\end{equation}
\end{thm}
\begin{pf}
Fix $\eta\in\real_{>0}^d$, and set $D=\dim\spn_\rat\{\eta_1,\ldots,\eta_d\}$.  As in Proposition~\ref{prop-positive-integer}, let $q\in\real_+^D$ and let $\alpha=(\alpha_1,\ldots,\alpha_d)$ be a sequence of multi-indices $\alpha_j\in\integer_+^D$ such that 
\[
\eta_j=\langle\alpha_j,q\rangle\qquad(1\leq j\leq d).
\]
Define $\sigma_\alpha$ as in (\ref{tau-alpha}), so that 
\[
\ell_\eta=\sigma_\alpha\circ\ell_q.
\]
Then $\tilde{f}\circ\sigma_\alpha\in\mathcal{S}_D$ by Proposition~\ref{prop-monomial-substitution}.  Using the notation (\ref{Lambda-notation}), it follows that 
\[
\Lambda \tilde{f}\circ\sigma_\alpha(1)=\sum_{j=1}^\infty-\log(1-|r_j|^2),
\]
by Theorem~\ref{thm-boyd} and the fact that $r$ and $\widetilde{r}$ have the same non-zero entries, where $\widetilde{r}$ is the sequence of Schur parameters of $\tilde{f}\circ\sigma_\alpha$ (see the proof of Proposition~\ref{prop-monomial-substitution}).   

Now, since the entries of $q$ are linearly independent over the integers, the torus line $\ell_q$ is dense on $\torus^D$.  Therefore Lemma~\ref{lem-continuous-birkhoff}, with
\[
g(z)=-\log\left(1-\bigl|\tilde{f}\circ\sigma_\alpha(z)\bigr|^2\right)\qquad(z\in\torus^d),
\]
implies that 
\[
\lim_{L\rightarrow\infty}\frac{1}{2L}\int_{-L}^L-\log\left(1-\bigl|(\tilde{f}\circ\sigma_\alpha)\circ\ell_q(\omega)\bigr|^2\right)d\omega=\Lambda \tilde{f}\circ\sigma_\alpha(1),
\]
completing the proof.
\end{pf}

\section{A trace formula for the Schr\"{o}dinger operator  with singular potential\label{sec-trace}}

The Sturm-Liouville problem on $[0,b]$, parameterized by $0\neq\omega\in\real$, 
\begin{gather}
\label{impedance-form}
(\zeta u^\prime)^\prime+\omega^2\zeta u=0\\
\label{boundary-conditions}
\textstyle
\frac{1}{2}\bigl(u(0)+\frac{1}{i\omega}u^\prime(0)\bigr)=1,\qquad\frac{1}{2}\bigl(u(b)-\frac{1}{i\omega}u^\prime(b)\bigr)=0
\end{gather}
determines a \emph{reflection coefficient} defined as 
\begin{equation}\label{reflection-coefficient}
\textstyle
R(\omega)=\frac{1}{2}\bigl(u(0)-\frac{1}{i\omega}u^\prime(0)\bigr). 
\end{equation}
Equation (\ref{impedance-form}) is called the Schr\"odinger equation in impedance form.  Provided $\zeta$ is sufficiently smooth, it relates by a well-known change of variables to the standard Schr\"odinger equation.  In detail, set 
\begin{equation}\label{alpha}
\alpha=-\zeta^\prime/(2\zeta)=\textstyle-\frac{1}{2}(\log\zeta)^\prime,
\end{equation} 
introduce a new dependent variable and coefficient 
\begin{equation}\label{y-q}
y=\zeta^{1/2}u,\qquad q=(\zeta^{1/2})^{\prime\prime}/\zeta^{1/2}=\alpha^2-\alpha^\prime, 
\end{equation}
and express (\ref{impedance-form}) in terms of $y$ to yield
\begin{equation}\label{schrodinger}
-y^{\prime\prime}+qy=\omega^2y.
\end{equation}
The potential $q$ extends to a function on $\real$ by setting $q(x)=0$ if $x\not\in[0,b]$, whereby the Schr\"odinger equation (\ref{schrodinger}) determines a scattering matrix defined in terms of Jost solutions (see \cite{DeTr:1979} for details).  The reflection coefficient (\ref{reflection-coefficient}) turns out to be the upper-right entry of this scattering matrix---and is hence the reflection coefficient in the standard sense for equation (\ref{schrodinger}). (See \cite[\S1.1.3]{Gi:Pre2021} for details.) 

The potential $q$ defined by (\ref{y-q}) has no ground states \cite[Thm.~1.1, p.~3092]{KaPeShTo:2005}.  It therefore follows from the first trace formula for the Schr\"odinger operator 
\cite[fla.~(1.1)]{HrMy:2021} that if $q\in L^1(\real)$,  
\begin{equation}\label{trace-1}
-\int_{\real}\log(1-|R(\omega)|^2)\,d\omega=\pi\int_\real q.
\end{equation}
On the other hand, if $q\not\in L^1(\real)$, 
then the trace formula (\ref{trace-1}) no longer applies.

The present section is concerned with the trace formula when the impedance function $\zeta$ is the restriction to $[0,b]$ of a strictly positive step function of the form
\begin{equation}\label{zeta-form}
\tilde\zeta=a_0\chi_{(-\infty,y_1)}+\left(\sum_{j=1}^{d-1}a_j\chi_{[y_j,y_{j+1})}\right)+a_d\chi_{[y_d,\infty)},\quad\mbox{ where }\quad 0<y_1<\cdots<y_d<b. 
\end{equation}
In this case, $\alpha$ is a purely distributional sum of Dirac functions, so the standard Schr\"odinger equation (\ref{schrodinger}) ceases to have a meaningful interpretation (even in the sense of distributions), and the formula (\ref{trace-1}) breaks down altogether.  
However, the impedance form of the Schr\"odinger equation 
remains well defined in the singular case (\ref{zeta-form}), as does the reflection coefficient (\ref{reflection-coefficient}).  And the almost periodic Szeg\H{o}-Verblunsky theorem of the previous section
supplies a singular analogue to (\ref{trace-1}), as follows.  
\begin{thm}\label{thm-singular-trace}
Let $R$ denote the reflection coefficient (\ref{reflection-coefficient}) for the impedance form (\ref{impedance-form}) of the Schr\"odinger equation in the case where $\zeta$ is positive of the form (\ref{zeta-form}), and set 
\begin{equation}\label{r-delta}
r_j=\frac{a_{j-1}-a_j}{a_{j-1}+a_j}\qquad(1\leq j\leq d).
\end{equation}
Then
\begin{equation}\label{trace-2}
\lim_{L\rightarrow\infty}\frac{1}{2L}\int_{-L}^L\log(1-|R(\omega)|^2)\,d\omega=\sum_{j=1}^d\log(1-|r_j|^2). 
\end{equation}
\end{thm}
\begin{pf}
Set $y_0=0$,
$
\eta_j=y_j-y_{j-1}$ $(1\leq j\leq d),
$
and write $\eta=(\eta_1,\ldots,\eta_d)$.  For piecewise constant $\zeta$, the reflection coefficient for (\ref{impedance-form}) can be computed explicitly as
$
R=f\circ\ell_\eta,
$
where, referring to (\ref{multiplication},\ref{composition}), 
\[
f=\mathfrak{m}_1\mathfrak{g}_{r_1}\cdots\mathfrak{m}_d\mathfrak{g}_{r_d}\mathbf{c}_0. 
\]
(See \cite[\S2]{Gi:JFAA2017} or \cite[\S3.5.2]{FoGaPaSo:2007} for details.)  Note that the Schur parameters for $f$ are $r_j=0$ for all $j>d$. 
Moreover, $f$ extends to a holomorphic function $\tilde{f}:\overline{\ddd}^d\rightarrow\ddd$, by (\ref{fn-formula-2}) and Proposition~\ref{prop-zero-free}. The desired formula (\ref{trace-2}) then follows from Theorem~\ref{thm-almost-periodic-main}. 
\end{pf}


\end{document}